%% file: 1238.tex
\documentclass[a4paper]{amsart} 
\usepackage[ascii]{inputenc} 

\usepackage{amssymb}
\usepackage{mathrsfs}
\usepackage[hidelinks]{hyperref}

\usepackage{xcolor}

\newcommand{\magenta}[1]{\textcolor{magenta}{#1}}

\usepackage[shortlabels]{enumitem}
\setlist[enumerate,1]{label={(\Alph*)}}
\setlist[enumerate,2]{label={(\alph*)}}
\setlist[enumerate,3]{label={$\bullet_{\arabic*}$}}

\newenvironment{PROOF}[2][\proofname.]
   {\begin{proof}[#1]\renewcommand{\qedsymbol}{\textsquare$_{\mathrm{#2}}$}}
   {\end{proof}}

\newtheorem{theorem}{Theorem}[section] 
\newtheorem{claim}[theorem]{Claim}

\newtheorem{observation}[theorem]{Observation}

\theoremstyle{definition}
\newtheorem{definition}[theorem]{Definition}

\newtheorem{example}[theorem]{Example}

\newtheorem{discussion}[theorem]{Discussion}
\newtheorem{convention}[theorem]{Convention}

\theoremstyle{remark}
\newtheorem{remark}[theorem]{Remark}
\newtheorem{question}[theorem]{Question}
\newtheorem{notation}[theorem]{Notation}

\input{mgrimes.tex}

\begin{document}
\makeatletter\def\shfiuwefootnote{\gdef\@thefnmark{}\@footnotetext}\makeatother\shfiuwefootnote{Version 2023-04-24\_2. See \url{https://shelah.logic.at/papers/1238/} for possible updates.}

\title {AEC for strictly stable \\
 Sh1238}
\author {Saharon Shelah}
\address{Einstein Institute of Mathematics\\
Edmond J. Safra Campus, Givat Ram\\
The Hebrew University of Jerusalem\\
Jerusalem, 91904, Israel\\
 and \\
 Department of Mathematics\\
 Hill Center - Busch Campus \\ 
 Rutgers, The State University of New Jersey \\
 110 Frelinghuysen Road \\
 Piscataway, NJ 08854-8019 USA}
\email{shelah@math.huji.ac.il}
\urladdr{http://shelah.logic.at}
\thanks{For earlier versions, the author thanks Alice Leonhardt for the beautiful typing up to 2019.
The author would like to thank the Israel Science Foundation for
partial support of this research (Grant No. 242/03). 
First typed 2002-Sept-10 as \S4-7 of \cite{Sh:839}; \S8-11 of the original paper can be found as \cite{Sh:1239}.\\
For later versions, the author would like to thank the NSF and BSF for partially supporting this research --- NSF-BSF 2021: grant with M. Malliaris, NSF 2051825, BSF 3013005232 (2021/10-2026/09). The author is also grateful to an individual who wishes to remain anonymous for generously funding typing services, and thanks Matt Grimes for the careful and beautiful typing.\\
References like [Sh:950, Th0.2=Ly5] mean that the internal label of Th0.2 is y5 in Sh:950.
The reader should note that the version in my website is usually more up-to-date than the one in arXiv.}

\makeatletter
\@namedef{subjclassname@2020}{\textup{2020} Mathematics Subject Classification}
\makeatother
\subjclass[2020]{Primary 03C48; Secondary 03C45}

\keywords {model theory, classification theory, stability, AEC, orthogonality, weight, main gap}

\date{April 23, 2023}

\begin{abstract}
Good frames were suggested in \cite{Sh:h} as 
the (bare-bones) parallel, in the context of AECs, to superstable 
(among elementary classes).  Here we consider $(\mu,\lambda,\kappa)$-frames as
candidates for being (in the context of AECs) the correct parallel 
to the class of $|T|^+$-saturated 
models of a strictly stable theory (among elementary classes).  
One thing we lose compared to the superstable case is that going up
by induction on cardinals is problematic (for stages of small
cofinality).  But this arises only when we try to lift such classes 
to higher cardinals. Also, we may use, as a replacement, the existence of prime models over unions of increasing chains. For this
context we investigate the dimension.
\end{abstract}

\maketitle
\numberwithin{equation}{section}
\setcounter{section}{-1}

\centerline {Annotated Content}
\bigskip


\noindent
\S0 \quad Introduction \hfill pg.\pageref{4}
\bigskip

\noindent
\S1 \quad Axiomatize AEC without full continuity (label f) \hfill pg.\pageref{5}
\sn
\begin{enumerate}
    \item[]  [Smooth out: generalize \cite[\S1]{Sh:600}.]
\end{enumerate}

\sn
\S(1A) \quad AEC \hfill pg.\pageref{5A}

\noindent
\S(1B) \quad Basic Notions (\ref{f31} on) \hfill pg.\pageref{5B}

\noindent
\S(1C) \quad Lifting such classes to higher cardinals (\ref{f46} on) \hfill pg.\pageref{5C}
\bigskip

\noindent
\S2 \quad PR frames (label g) \hfill pg.\pageref{6}
\mn
\bigskip

\noindent
\S3 \quad Thoughts on the main gap (label h) \hfill pg.\pageref{7}

\newpage


\section{Introduction} \label{4}

In this part we try to deal with classes like ``$\aleph_1$-saturated
models of a first order theory $T$, and even strictly stable ones" rather
than of ``a model of $T$," in the AEC framework.  
The parallel problem for ``model of $T$, even superstable one"
is the subject of \cite{Sh:h}.

\medskip
Now, some constructions go well by induction on cardinality (say, by
dealing with {a} $(\lambda,\clP(n))$-system of models) but not all.
E.g., starting with $\aleph_0$ we may consider $\lambda > \aleph_0$, so we
can find $F:[\lambda]^{\aleph_0} \to \lambda$ such that there
is no infinite decreasing sequence of $F$-closed subsets of $\lambda$ 
{such that} $u \in [\lambda]^{< \aleph_0} \Rightarrow F(u) = \varnothing$, but maybe such that $u \in
[\lambda]^{\le \aleph_0} \Rightarrow |\cl_F(u)| \le \aleph_0$.
Let $\LL u_\alpha:\alpha < \alpha_* \RR$ list $\{\cl_F(u):u
\in [\lambda]^{\le \aleph_0}\}$ such that $\cl_F(u_\alpha)
\subseteq \cl_F(u_\beta) \Rightarrow \alpha \le \beta$. We try to
choose $M_{u_\alpha}$ by induction on $\alpha$.

Another approach is to consider strictly stable theories with `nice enough' type theory, like superstable. See \cite{Sh:1133} on so-called `flat' first order theories.
\smallskip

\centerline {$* \qquad * \qquad *$}

\medskip
This was the middle part of \cite{Sh:839} and was divided by editor request; the third part is \cite{Sh:1239}. The original full paper has existed (and to some extent, has circulated) since 2002.

\begin{notation}\label{z2}
Let $\lambda^{<\kappa} \defeq \sum\{\lambda^\sigma : \sigma < \kappa\}$; in subscripts we may use $\lambda[<\kappa]$. 
\end{notation}

We shall use {the following} freely.

\begin{claim}\label{z5}
If $\lambda = \lambda^{<\kappa}$ then $\chi \geq \lambda \Rightarrow (\chi^{<\kappa})^{<\kappa} = \chi^{<\kappa}$.
\end{claim}

\newpage

\section{Axiomatizing AEC without full continuity} \label{5}

\subsection{DAEC} \label{5A} 
\bigskip

Classes like ``the $\aleph_1$-saturated models of a first order $T$ which
is not superstable", do not fall under AEC --- still, they are close,
and below we suggest a framework for them.  So for increasing sequences of
short length the union is not necessarily in the class, but we have weaker demands. In the main case, as `compensation,' we have {that} prime models exist; 
in particular, over short increasing chains of models.
\smallskip


We shall lift a $(\mu,\lambda,\kappa)$-AEC to
$(\infty,\lambda,\kappa)$-AEC (see below), so actually 
$\gk_\lambda$ will suffice. But for our main objects -- good frames -- this is
more complicated, as their properties (e.g., the amalgamation property)
are not necessarily preserved by the lifting.

This section generalizes \cite[\S1]{Sh:600}; in some cases the
differences are minor, whereas sometimes the differences are the whole point.

\begin{convention}\label{f0}
In this section, if not said otherwise, $\gk$ will denote a 1-DAEC (i.e. a directed AEC; see Definition
\ref{f2}). We may write DAEC (the \emph{D} stands for directed).

\end{convention}

\begin{definition}\label{f2}  
Assume $\lambda < \mu$, $\lambda^{<\kappa} = \lambda$ (for notational simplicity),\\ 
$\alpha < \mu \Rightarrow |\alpha|^{< \kappa} < \mu$, and $\kappa$ is regular.

We say that $\gk$ is a $(\mu,\lambda,\kappa)$-1-DAEC when $\boxplus$ and all the axioms below hold.

(We may omit or add the `1' and/or the `$(\mu,\lambda,\kappa)$' by $\boxplus$(a) below; similarly in similar definitions.
Instead {of} $\mu = \mu^+_1$ we may write $\le \mu_1$.)
We write pre-DAEC or 0-DAEC when we omit \textbf{Ax.III}(b), \textbf{IV}(b). 
\mn
\begin{enumerate}
    \item[$\boxplus$] (=\textbf{Ax.O}) $\gk$ consists of the objects in clauses (a)-(d), having the properties listed in (e)-(g).
    \begin{enumerate}
        \item The cardinals $\mu = \mu_\gk = \mu(\gk)$, $\lambda = \lambda_\gk = \lambda(\gk)$ 
        and $\kappa = \kappa_\gk = \kappa(\gk)$, satisfying 
        $\mu > \lambda = \lambda^{< \kappa} \ge \kappa = \cf(\kappa)$ 
        and $\alpha < \mu \Rightarrow |\alpha|^{< \kappa} < \mu$ (but possibly $\mu = \infty$).
\sn
        \item $\tau_\gk$, a vocabulary with each predicate and function symbol of arity $\le \lambda$.

        \item $K$ a class of $\tau$-models. 
\sn
        \item A two-place relation $\le_\gk$ on $K$.
\sn
        \item[(e)]  If $M_1 \cong M_2$ then $M_1 \in K \Leftrightarrow M_2 \in K$.
\sn
        \item[(f)]  if $(N_1,M_1) \cong (N_2,M_2)$ \underline{then} $M_1 \le_\gk N_1 \Rightarrow M_2 \le_\gk N_2$.
\sn
        \item[(g)]  Every $M \in K$ has cardinality $\lambda \leq \|M\| < \mu$.
    \end{enumerate}   
\end{enumerate}

\begin{enumerate}[align=left, labelwidth=1.25cm, leftmargin=\dimexpr\labelwidth+\labelsep,itemindent=!]
    \item[\textbf{Ax.I}(a)]  $M \le_\gk N \Rightarrow M \subseteq N$
\sn
    \item[\textbf{Ax.II}(a)]  $\le_\gk$ is a partial order.
\sn
    \item[\textbf{Ax.III}]   Assume that $\LL M_i : i < \delta \RR$ is a 
    $\le_\gk$-increasing sequence and\\ $\big\| \bigcup\{M_i : i < \delta\} \big\| < \mu$. 
    \underline{Then}:
\sn
    \begin{enumerate}
        \item   \textbf{Existence of unions} 
        
        If $\cf(\delta) \ge \kappa$ then there is $M \in K$ such that 
        $i < \delta \Rightarrow M_i \le_\gk M$ and $|M| = \bigcup\{|M_i|:i < \delta\}$
        (but not necessarily $M = \bigcup\limits_{i < \delta} M_i$).
\sn
        \item  \textbf{Existence of limits} 
        
        There is $M \in K$ such that $i < \delta \Rightarrow M_i \le_\gk M$.
    \end{enumerate}
    
    \item[\textbf{Ax.IV}(a)]  \textbf{Weak uniqueness of limit} (= weak smoothness)

    \noindent
    For $\LL M_i : i < \delta\RR$ as above, 
\sn
    \begin{enumerate}
        \item  If $\cf(\delta) \ge \kappa$, $M$ is as in \textbf{Ax.III}(a), and 
        $i < \delta \Rightarrow M_i \le_\gk N$, then $M \le_\gk N$. (This implies the uniqueness of $M$.)
\sn
        \item  If $N_\ell \in K$ and $i < \delta \Rightarrow M_i \le_\gk N_\ell$ 
        for $\ell =1,2$ \underline{then} there are $N \in K$ and $f_1,f_2$ such 
        that $f_\ell$ is a $\le_\gk$-embedding of $N_\ell$ into $N$ for $\ell=1,2$ 
        and $i < \delta \Rightarrow f_1 \rest M_i = f_2 \rest M_i$.
    \end{enumerate}

    \item[\textbf{Ax.V}]   If $N_\ell \le_\gk M$ for $\ell=1,2$ and 
    $N_1 \subseteq N_2$ \underline{then} $N_1 \le_\gk N_2$.
\sn
    \item[\textbf{Ax.VI}]   \textbf{L.S.T. property}
    
    \noindent If $A \subseteq M \in K$ and $|A| \le \lambda$ \underline{then} there is $M \le_\gk N$ of cardinality $\lambda$ such that $A \subseteq M$.
\end{enumerate}
\end{definition}

\begin{remark}\label{f3}
There are some more axioms listed in \ref{f5}(5), but
we shall mention them in any claim in which they are used so no need
to memorize. Note that \ref{f5}(1)-(4) assumes some of them.
\end{remark}

\begin{definition}\label{f5}
1) We say $\gk$ is a 4-DAEC or DAEC$^+$ \underline{when} it is a
$(\lambda,\mu,\kappa)$-1-DAEC and satisfies \textbf{Ax.III}(d), \textbf{Ax.IV}(e) below.

\noindent
2) We say $\gk$ is a 2-DAEC or DAEC$^\pm$ \underline{when} it is a
$(\lambda,\mu,\kappa)$-0-DAEC and \textbf{Ax.III}(d), \textbf{Ax.IV}(d) below hold.

\noindent 
3) We say $\gk$ is 5-DAEC \underline{when} it is 1-DAEC and \textbf{Ax.III}(d),(f) holds.

\noindent
4) We say $\gk$ is 6-DAEC \underline{when} it is a 1-DAEC and
   \textbf{Ax.III}(d),(f) + \textbf{Ax.IV}(f).

\noindent
5) Concerning Definition \ref{f2}, we consider the following axioms:

\mn
\begin{enumerate}
    \item[\textbf{Ax.III}] 
    \begin{enumerate}
        \item[(c)] If $I$ is $\kappa$-directed and $\olsi M = \LL M_s : s \in I\RR$ 
        is $\le_\gk$-increasing (that is, $s \le_I t \Rightarrow M_s \leq_\gk M_s$), and 
        $\sum\{\|M_s\| : s \in I\} < \mu$ \underline{then} $\olsi M$ has a  $\le_\gk$-upper bound $M$ (i.e. $s \in I \Rightarrow M_s \le_\gk M$).
\sn
        \item[(d)] \textbf{Union of directed systems} 
        
        If $I$ is $\kappa$-directed, 
        $|I| < \mu$, $\LL M_t : t \in I \RR$ is $\le_\gk$-increasing, and\\ 
        $\big\| \bigcup \{M_t : t \in I\} \big\| < \mu$ \underline{then} there is one 
        and only one $M$ with universe $\bigcup\{|M_t|:t \in I\}$ such that 
        $M_s \le_\gk M$ for every $s \in I$. 
        (We call it the $\le_\gk$-union of $\LL M_t : t \in I\RR$.)
\sn
        \item[(e)] Like \textbf{Ax.III}(c), but $I$ is just directed.
\sn
        \item[(f)] If $\olsi M = \LL M_i : i < \delta\RR$ is $\le_\gk$-increasing, 
        $\cf(\delta) < \kappa$, and 
        $$\big|\textstyle\bigcup\{M_i : i < \delta\} \big| < \mu$$
        \underline{then} there is $M$ which is $\le_\gk$-prime over $\olsi M$; i.e.
        \begin{itemize}
            \item If $N \in K_\gk$ and $i < \delta \Rightarrow M_i \le_\gk N$ 
            \underline{then} there is a $\le_\gk$-embedding of $M$ into 
            itself over $\bigcup\{|M_i| : i < \delta\}$.
        \end{itemize}
    \end{enumerate}

    \item[\textbf{Ax.IV}]
    \begin{enumerate}
        \item[(c)] If $I$ is $\kappa$-directed and $\olsi M = \LL M_s : s \in I\RR$ 
        is $\le_\gk$-increasing and $N_1,N_2$ are $\le_\gk$-upper bounds of $\olsi M$ \underline{then} for some $(N'_2,f)$ we have $N_2 \le_\gk N'_2$ and $f$ is a $\le_\gk$-embedding of $N_1$ into $N_2$ which is the identity on $M_s$ for every $s \in I$.
        (This is a weak form of uniqueness.)
\sn
        \item[(d)] If $I$ is a $\kappa$-directed partial order, $\olsi M = \LL M_s : s \in I\RR$ 
        is $\le_\gk$-increasing, $s \in I \Rightarrow M_s \le_\gk M$ and 
        $|M| = \bigcup\{|M_s|:s \in I\}$, \underline{then} 
        $\bigwedge\limits_{s} M_s \le_\gk N \Rightarrow M \le_\gk N$.
\sn
        \item[(e)] Like \textbf{Ax.IV}(c), but $I$ is just directed.
\sn
        \item[(f)] If $I$ is directed and $\olsi M = \LL M_s : s \in I\RR$ is $\le_\gk$-increasing \underline{then} there is $M$ which is a $\le_\gk$-prime over $\olsi M$, defined as in \textbf{Ax.III}(f).
    \end{enumerate}
\end{enumerate}
\end{definition}

\begin{claim}\label{f6}
Assume\footnote{By \ref{f0}, it is not necessary to say this.} $\gk$ is a DAEC.

\noindent
1) \textbf{Ax.III}(d) implies \textbf{Ax.III}(c). Also, \textbf{Ax.III}(c) implies \textbf{Ax.III}(a).

\noindent
2) \textbf{Ax.III}(e) implies \textbf{Ax.III}(c) and also \textbf{Ax.III}(b).
Also, \textbf{Ax.III}(b) implies \textbf{Ax.III}(a).

\noindent
3) \textbf{Ax.IV}(d) implies \textbf{Ax.IV}(a).

\noindent
4) \textbf{Ax.IV}(e) implies \textbf{Ax.IV}(c) and also \textbf{Ax.III}(b).

\noindent
5) In all the axioms in Definition \ref{f5} it is necessary that
$\big|\bigcup\{M_s:s \in I\}\big| < \mu_\gk$.

\noindent
6) \textbf{Ax.IV}(b) implies that $\gk$ has amalgamation.
\end{claim}

\begin{PROOF}{\ref{f6}}
Easy.
\end{PROOF}

\begin{example}\label{f7}
\textbf{The first order case}

Let $T$ be a stable complete first order theory, and $\kappa = \kappa_r(T) \in \big[ \aleph_1, |T|^+ \big]$. (Equivalently, $\kappa$ is the minimal regular cardinal such that $\lambda = \lambda^{<\kappa} \geq 2^{|T|} \Rightarrow T \text{ stable in } \lambda$.) We shall define $\gk = \gk_T$:
\begin{enumerate}
    \item[$\boxplus$] 
    \begin{enumerate}
        \item $K = K_\gk$ is the class of $\kappa$-saturated models of $T$ (equivalently, $\bfF_\kappa^a$-saturated).

        \item ${\leq} = {\leq_\gk}$ means ``is an elementary submodel of.''
    \end{enumerate}
\end{enumerate}
\end{example}

\begin{example}\label{f9}
\textbf{Existentially closed}

Let $T$ be a universal first order theory with the JEP, for transparency.\footnote{Otherwise the class of existentially closed models of $T$ is divided into $\leq 2^{|T|}$ subclasses, each of them of this form.} We shall define $\gk$:
\begin{enumerate}
    \item[$\boxplus$] 
    \begin{enumerate}
        \item $K = K_\gk$ is the class of existentially closed models of $T$.

        \item ${\leq} = {\leq_\gk}$ means ``is a submodel of.''
    \end{enumerate}
\end{enumerate}
See, in \cite{Sh:54}, what was called the ``kind III context''; recall that kind II was for such $T$ with amalgamation and JEP. Much more was done by Hrushovski.
\end{example}

\begin{example}\label{f12}
\textbf{Metric Spaces}

\begin{enumerate}[1)]
    \item 
    \begin{enumerate}
        \item We say that $\tau$ is a \emph{metric vocabulary} if if has the distinguished 2-place predicates $R_q$ ($q$ a positive rational) and $\nmet$:\footnote{nmet stands for `non-metric.'} let
        $$\nmet(\tau) = \tau \setminus \{ R_q : q \in \bbQ^+\}$$ 
        and $\tau$ is finitary. That is, each predicate and function symbol has finitely many places.

        \item $M$ is a \emph{metric model} when its vocabulary $\tau_M$ is a metric vocabulary and there is a metric $\bfd_M(-,-)$ on $M$ such that 
        \begin{enumerate}
            \item $\bfd_M(a,b) = \inf\{q \in \bbQ^+ : (a,b) \in R_q^+\}$

            \item For any predicate $R \in \tau$, $R^M$ is closed.

            \item For any function symbol $F$, $F^M$ is a continuous function.

            \item $M$ is complete as a metric space.
        \end{enumerate}

        \item Without clause $\bullet_4$, we say $M$ is an \emph{almost metric} model.

        \item We say that the metric models $M_1,M_2$ are \emph{topologically isomorphic} when there is a $\pi$ such that
        \begin{enumerate}
            \item $\pi$ is an isomorphism from $M_1 \rest \nmet(\tau)$ onto $M_2 \rest \nmet(\tau)$.

            \item $\dist_\pi(M_1,M_2) \defeq$ 
            $$\sup \!\bigg\{ \frac{\bfd_{M_2}\big(\pi(a),\pi(b)\big)}{ \bfd_{M_1}(a,b)}, \frac{\bfd_{M_1}(a,b)}{\bfd_{M_2}\big(\pi(a),\pi(b)\big)} : a \neq b \in M_1 \bigg\}$$ 
            is finite.
        \end{enumerate}
        (Note that this is the meaning of isomorphism for Banach space theorists; what we call isomorphism they would call isometry.)
    \end{enumerate}

    \item 
    \begin{enumerate}
        \item We say $\gk$ is a \emph{metric} AEC (or MAEC) \underline{when}:
        \begin{enumerate}
            \item $\tau_\gk$ is a metric vocabulary.

            \item $\gk$ is a DAEC with $\mu_\gk = \infty$, $\kappa = \aleph_1$, and $\lambda = \lambda^{\aleph_0}$ (and for convenience $|\tau_\gk| \leq \lambda$).

            \item Each $M \in K_\gk$ is a metric model.

            \item If $I$ is a directed partial order and $\olsi M = \LL M_s : s \in I\RR$ is $\leq_\gk$-increasing \underline{then} the completion $M$ of $\bigcup\{ M_s : s \in I\}$, naturally defined, is a $\leq_\gk$-l.u.b. of $\olsi M$.
        \end{enumerate}

        \item We say $\gk$ is an \emph{almost metric} AEC when we omit the completeness demand in (1)(b), and add
        \begin{itemize}
            \item If $N$ is the completion of $M \in K_\gk$ (so necessarily $N \in K_\gk$, $M \leq_\gk {N}$) then $M \subseteq M' \subseteq N \Rightarrow M \leq_\gk M' \leq_\gk N$. 
        \end{itemize}
    \end{enumerate}

    \item
    \begin{enumerate}
        \item If $\gk$ is a metric AEC \underline{then} all the axioms in Definitions \ref{f2}, \ref{f5} hold.        

        \item If $\gk$ is an almost metric AEC \underline{then} 
        $$\comp(\gk) \defeq \gk \rest \big\{ M \in K_\gk : \big(|M|,\bfd_M \big) \text{ is complete}\big\}$$
        is a metric AEC; also, $\gk$ is an AEC. In this case, ``the completion of $M \in K_\gk$'' is naturally defined.

        \item \textbf{The representation theorem}.

        If $\gk$ is a metric AEC \underline{then} for some $\tau_1,T_1,\Gamma$, we have:
        \begin{enumerate}
            \item $\tau_1 \supseteq \tau_\gk$ and $|\tau_1| \leq \lambda_\gk$.

            \item $T_1$ is a universal f.o. theory in $\bbL(\tau_1)$.

            \item $\Gamma$ is a set of $\bbL(\tau_1)$-types consisting of formulas (so they are $m$-types for some $m$).

            \item Every $M \in \EC(T_1,\Gamma)$ is a weak metric model.

            \item $K_\gk = \{M : M \text{ is the completion of $M_1 \rest \tau_\gk$ for some } M_1 \in \EC(T_1,\Gamma)\}$

            \item $\leq_\gk$ is defined as
            \begin{align*}
               \big\{(M,N) : &\ \text{there are $M_1 \subseteq N_1$ from $\EC(T_1,\Gamma)$ such that}\\ 
               &\ M \subseteq N \text{ are the completions of } M_1 \rest \tau_\gk,\ N_1 \rest \tau_\gk \text{ resp.}\big\} 
            \end{align*}
        \end{enumerate} 
    \end{enumerate}
    [Why is this true? As in the AEC case.]
\end{enumerate}
Regarding metric model theory and topological model theory, the field was started by Chang and Keisler in \cite{ChKe66}; for an introduction see a recent survey \cite{Kei20}.
\end{example}

\begin{definition}\label{f15}
We say $\LL M_i : i < \alpha \RR$ is 
$\le_\gk$-increasing $(\ge \kappa)$-continuous
\underline{when} it is $\le_\gk$-increasing and $\delta < \alpha$ and 
$\cf(\delta) \ge \kappa \Rightarrow |M_\delta| = \bigcup\{|M_j| : j < \delta\}$.
\end{definition}

\noindent
As an exercise we consider directed systems with mappings.

\begin{definition}\label{f16}
1) We say that $\olsi M = \LL M_t,h_{t,s} : s \le_I t\RR$ is 
a $\le_\gk$-directed system \underline{when}
\begin{enumerate}
    \item $I$ is a directed partial order.

    \item If $s \leq_I t$ then $h_{t,s}$ is an isomorphism from $M_s$ onto some $M' \leq_\gk M_t$.
 
    \item If $t_0 \le_I t_1 \le_I t_2$ then 
    $h_{t_2,t_0} = h_{t_2,t_1} \circ h_{t_1,t_0}$.
\end{enumerate}

\sn
1A) We say that $\olsi M = \LL M_t,h_{t,s}:s \le_I t \RR$ is a
$\le_\gk$-$\theta$-directed system when in addition $I$ is
$\theta$-directed.

\sn
2) We may omit $h_{t,s}$ when $s \le_I t \Rightarrow h_{t,s} = 
\id_{M_s}$ and write $\olsi M = \LL M_t : t \in I\RR$.

\sn
3) We say $(M,\bar h)$ is a $\le_\gk$-limit of $\olsi M$ when 
$\bar h = \LL h_s : s \in I\RR$, $h_s$ is a 
$\le_\gk$-embedding of $M_s$ into $M$, and 
$s \le_I t \Rightarrow h_s = h_t \circ h_{t,s}$.

\sn
4) We say $\olsi M = \LL M_\alpha:\alpha < \alpha^*\RR$ is
 $\le_\gk$-semi-continuous \underline{when}: (see \textbf{Ax.III}(f) in \ref{f5})
\mn
\begin{enumerate}
    \item   $\olsi M$ is $\le_\gk$-increasing.
\sn
    \item If $\alpha < \alpha^*$ has cofinality $\ge \kappa$ then 
    $M_\alpha = \bigcup\{M_\beta : \beta < \alpha\}$.
\sn
    \item  If $\alpha < \alpha^*$ has cofinality $<\kappa$ then 
    $M_\delta$ is $\le_\gk$-prime over $\olsi M \rest \alpha$.
\end{enumerate}
\end{definition}

\begin{observation}\label{f17}
$[\gk$ is a DAEC.$]$

\sn
1) If $\olsi M = \LL M_t,h_{t,s}:s \le_I t \RR$ is a 
$\le_\gk$-directed system, \underline{then} we can
find a $\le_\gk$-directed system $\LL M'_t:t \in I\RR$
(so $s \le_I t \Rightarrow M'_s \le_\gk M'_t$) and $\bar g =
\LL g_t:t \in I\RR$ such that:
\mn
\begin{enumerate}
    \item[$(a)$]   $g_t$ is an isomorphism from $M_t$ onto $M'_t$.
\sn
    \item[$(b)$]   If $s \le_I t$ then $g_s = g_t \circ h_{t,s}$.
\end{enumerate}
\mn
2) So in the axioms \textbf{III}(a),(b),\textbf{IV}(a) from Definition \ref{f2} as
well as those of \ref{f5} we can use $\le_\gk$-directed system
$\LL M_s,h_{t,s}:s \le_I t \RR$ with $I$ as there.

\sn
3) If $\gk$ is an ess-$(\mu,\lambda)$-AEC (see \cite[\S1]{Sh:839}) \underline{then} $\gk$
   is a $(\mu,\lambda,\aleph_0)$-DAEC and satisfies all the axioms
   from \ref{f5}.

\sn
4) If $(M,\bar h)$ is prime over $\olsi M = \LL
M_t,h_{t,s}:s \le_I t\RR$ and $\chi = \sum\{\|M_t\|:t \in I\}$
then $\|M\| \le \chi^{< \kappa}$.
\end{observation}

\begin{PROOF}{\ref{f17}}
Straightforward; e.g. we can use ``$\gk$ has $(\chi^{< \kappa})$-LST"
(i.e. Observation \ref{f19} below).
\end{PROOF}

\noindent
More serious is proving the LST theorem in our context (recall that in
the axioms, see \textbf{Ax.VI}, we demand it only down to $\lambda$).

\begin{claim}\label{f19}
1) $\gk$ is a $(\mu,\lambda,\kappa)$-2-DAEC --- see Definition \ref{f5}(2).

If $\lambda_\gk \le \chi = \chi^{< \kappa} < \mu_\gk$, $A \subseteq N \in \gk$, 
and $|A| \le \chi \le \|N\|$ \underline{then} there is $M \le_\gk N$ 
of cardinality $\chi$ such that $\|M\| = \chi$ and $A \subseteq M$.

\sn 
2) If $\gk$ satisfies \textbf{Ax.III}(e), \textbf{Ax.IV}(e) \underline{then} in part (1) we do not need the assumption $``\chi = \chi^{<\kappa}."$ 
\end{claim}

\begin{PROOF}{\ref{f19}}
1) As $\chi \leq \|N\|$,
\begin{enumerate}
    \item[$(*)_0$] Without loss of generality $|A| = \chi$.
\end{enumerate}

Let $\big\LL u_\alpha : \alpha < \alpha(*)\big\RR$ list 
$[A]^{< \kappa(\gk)}$ and let $I$ be the following partial order:
\mn
\begin{enumerate}
    \item[$(*)_1$] 
    \begin{enumerate}
        \item The set of elements is $\{\alpha < \chi : 
        \text{for no $\beta < \alpha$ do we have } u_\alpha \subseteq u_\beta\}$.
\sn
        \item $\alpha \le_I \beta$ iff $u_\alpha \subseteq u_\beta$ 
        (hence $\alpha \le \beta$).
    \end{enumerate}
\end{enumerate}
\mn
Easily
\mn
\begin{enumerate}
    \item[$(*)_2$]
    \begin{enumerate}
        \item $I$ is $\kappa$-directed.

        \item For every $\alpha < \alpha(*)$, for some $\beta < \alpha(*)$, 
        we have $u_\alpha \subseteq u_\beta \wedge \beta \in I$.

        \item $\bigcup\{u_\alpha:\alpha \in I\} = A$.
    \end{enumerate}
\end{enumerate}
\mn
Now we choose $M_\alpha$ by induction on $\alpha < \chi$ such that
\mn
\begin{enumerate}
    \item[$(*)_3$]   
    \begin{enumerate}
        \item $M_\alpha \le_\gk N$

        \item $\|M_\alpha\| = \lambda_\gk$

        \item $M_\alpha$ includes $\bigcup\{M_\beta : \beta <_I \alpha\} \cup u_\alpha$.
    \end{enumerate}  
\end{enumerate}
\mn
Note that 
$$\big|\{\beta \in I:\beta <_I \alpha\}\big| \le \big|\{u : u \subseteq u_\alpha\}\big| = 2^{|u_\alpha|} \le 2^{< \kappa(\gk)} \le \lambda_\gk,$$ 
and by the induction hypothesis $\beta < \alpha \Rightarrow \|M_\beta\| \le \lambda_\gk$. 
Recall $|u_\alpha| < \kappa(\gk) \le \lambda_\gk$ hence the set 
$\bigcup\{M_\beta : \beta < \alpha\} \cup u_\alpha$ is a subset of $N$ of cardinality 
$\le \lambda$, hence by \textbf{Ax.VI} there exists $M_\alpha$ as required.

Having chosen $\LL M_\alpha : \alpha \in I\RR$, clearly by \textbf{Ax.V} it is a
$\le_\gk$-increasing $(<\kappa)$-directed system; hence by \textbf{Ax.III}(d), 
$M = \bigcup\{M_\alpha : \alpha \in I\}$ is well defined 
with universe $\bigcup\{|M_\alpha| : \alpha \in I\}$ and by \textbf{Ax.IV}(d) we
have $M \le_\gk N$.

Clearly $\|M\| \le \sum\{\|M_\alpha\| : \alpha \in I\} \le |I| \cdot \lambda_\gk = \chi$, 
and by $(*)_2(c) + (*)_3(c)$ we have 
$$A \subseteq \bigcup\{u_\alpha : \alpha < \chi\} = \bigcup\{u_\alpha : \alpha \in I\}
\subseteq \bigcup\{|M_\alpha| : \alpha \in I\} = M$$ 
and so $M$ is as required.  

\sn
2) Similarly.
\end{PROOF}

\begin{notation}\label{f22}
1) For $\chi \in [\lambda_\gk,\mu_\gk)$ let 
$K_\chi = K^\gk_\chi = \{M \in K : \|M\| = \chi\}$ and
$K_{< \chi} = \bigcup\limits_{\mu < \chi} K_\mu$. 

\sn
2) $\gk_\chi = (K_\chi,\le_\gk \rest K_\chi)$.

\sn
3) If $\lambda_\gk \le \lambda_1 < \mu_1 \le \mu_\gk$, $\lambda_1 =
\lambda^{<\kappa}_1$, and $(\forall \alpha < \mu_1)\big[|\alpha|^{< \kappa} < \mu_1\big]$ 
\underline{then} we define $K_{[\lambda_1,\mu_1]} = K^\gk_{[\lambda_1,\mu_1]}$ 
and $\gk_1 = \gk_{[\lambda_1,\mu_1]}$ similarly: 
\begin{enumerate}
    \item $K_{\gk_1} = \big\{M \in K_\gk : \|M\| \in [\lambda_1,\mu_1)\big\}$

    \item ${\le_{\gk_1}} = {\le_\gk \rest K_{\gk_1}}$

    \item $\lambda_{\gk_1} = \lambda_1$, $\mu_{\gk_1} = \mu_1$, $\kappa_{\gk_1} = \kappa_\gk$.
\end{enumerate} 
4) Let $\gk_{[\lambda_1,\mu_1]} \defeq \gk_{[\lambda_1,\mu_1^+)}$.
\end{notation}

\begin{definition}\label{f24}
The embedding $f:N \to M$ is a $\gk$-embedding or a 
$\le_\gk$-embedding \underline{when} its range is the universe of a model 
$N' \le_\gk M$ (so $f:N \to N'$ is an isomorphism, hence it is onto).
\end{definition}

\begin{claim}\label{f27}
$[\gk$ is a 2-DAEC.$]$

\sn
1) For every $N \in K$ there is a $\kappa_\gk$-directed partial order $I$ 
of cardinality $\le \|N\| < \kappa$ and $\olsi M = \LL M_t : t \in I \RR$ such that 
$t \in I \Rightarrow M_t \le_\gk N$, $\|M_t\| \le \LST(\gk) = \lambda_\gk$, 
$$I \models ``s < t" \Rightarrow M_s \le_\gk M_t,$$
and $N = \bigcup\limits_{t \in I} M_t$. 

\sn
1A) If $\gk$ satisfies \textbf{Ax.III}(e), \textbf{Ax.IV}(e) \underline{then} in part (1) we can add $|I| \leq \|M\|$.

\sn
2) For every $N_1 \le_\gk N_2$ we can find $\LL M^\ell_t : t \in I^* \RR$ 
as in part (1) for $N_\ell$ such that $I_1 \subseteq I_2$ 
and $t \in I_1 \Rightarrow M^2_t = M^1_t$. 
\end{claim}

\begin{PROOF}{\ref{f27}}
1), 1A) As in the proof of \ref{f19}.

\sn
2) Similarly.
\end{PROOF}

\begin{claim}\label{f28}
Assume $\lambda_\gk \le \lambda_1 = \lambda^{< \kappa}_1 < \mu_1 \le \mu_\gk$ 
and $(\forall\alpha < \mu_1)\big[|\alpha|^{< \kappa} < \mu_1\big]$.

\sn
1) \underline{Then} $\gk^*_1 \defeq \gk_{[\lambda_1,\mu_1)}$ as defined in
\ref{f22}(3) is a $(\lambda_1,\mu_1,\kappa_\gk)$-DAEC.

\sn
2) For each of the following axioms, if $\gk$ satisfies it then so
does $\gk_1$: \textbf{Ax.III}(d),(e), \textbf{Ax.IV}(c),(d),(e).

\sn
3) In part (2), its conclusion also applies to \textbf{Ax.III}(f), \textbf{Ax.IV}(f).


\end{claim}

\begin{PROOF}{\ref{f28}}
Read the definitions.
\end{PROOF}

\begin{claim}
\label{f29}
1) If $\gk$ satisfies \textbf{Ax.IV}(e) \underline{then} $\gk$ satisfies \textbf{Ax.III}(e)
provided that $\mu_\gk$ is regular or at least the relevant $I$ has
cardinality $< \cf(\mu_\gk)$. 

\noindent
2) If \textbf{Ax.III}(d),\textbf{IV}(d) hold, we can waive `$\mu_\gk$ is regular.'
\end{claim}

\begin{PROOF}{\ref{f29}}
Recall $\gk$ is a DAEC. We prove this by induction on $\theta = |I|$. 

Let $\chi = \lambda + \theta +\sum \big\{ \|M_s\| : s \in I \big\}$, which is in $[\lambda,\mu)$.

\mn
\textbf{Case 1:}  $I$ is finite.

So there is $t^* \in I$ such that $t \in I \Rightarrow t \le_I t^*$,
so this is trivial.

\mn
\textbf{Case 2:}  $I$ is countable.

So we can find a sequence $\LL t_n : n < \omega \RR$ such that $t_n \in I$, 
$t_n \le_I t_{n+1}$, and $s \in I \Rightarrow \bigvee\limits_{n < \omega} s \le_I t_n$.  
Now we can apply \textbf{Ax.III}(b) to $\LL M_{t_n} : n < \omega \RR$.

\mn
\textbf{Case 3:} $I$ uncountable.

First, we can find an increasing continuous sequence 
$\LL I_\alpha : \alpha < |I| \RR$ such that
$I_\alpha \subseteq I$ is directed of cardinality $\le |\alpha| + \aleph_0$ 
and $I_{|I|} = I = \bigcup\{I_\alpha : \alpha < |I|\}$.

Second, by the induction hypothesis for each $\alpha < |I|$ we 
choose $N_\alpha$ and $\bar h^\alpha = \LL h_{\alpha,t} : t \in
I_\alpha \RR$ such that:
\sn
\begin{enumerate}
    \item   $N_\alpha \in \gk^\gs_{\le \chi} \subseteq \gk$
\sn
    \item  $h_{\alpha,t}$ is a $\le_\gk$-embedding of $M_t$ into $N_\alpha$.
\sn
    \item   If $s <_I t$ are in $I_\alpha$ then 
    $h_{\alpha,s} = h_{\alpha,t} \circ h_{t,s}$.
\sn
    \item If $\beta < \alpha$ then $N_\beta \le_\gk N_\alpha$ and $t \in I_\beta \Rightarrow h_{\alpha,t} = h_{\beta,t}$.
\end{enumerate}
\sn
For $\alpha = 0$ use the induction hypothesis.

For $\alpha$ a limit ordinal, by \textbf{Ax.III}(a) there is $N_\alpha$ as
required; as $I_\alpha = \bigcup\limits_{\beta < \alpha} I_\beta$, there are no
new $h_t$-s. (Well, we have to check $\sum\{\|N_\beta\| : \beta < \alpha\} < \mu_\gk$, 
but as we assume $\mu_\gk$ is regular this holds.)

For $\alpha = \beta +1$, by the induction hypothesis there is
$(N'_\alpha,\bar g^\alpha)$ which is a limit of $\LL M_s,h_{t,s}:s
\le_{I_\alpha} t\RR$.  Now apply \textbf{Ax.IV}(e): well, apply the directed
system version with $\LL M_s,h_{t,s} : s \le_{I_\beta} t\RR$, $(N'_\alpha,\bar g_\alpha)$, 
$(N_\beta,\LL h_s : s \in I_\beta\RR)$ here standing for $\olsi M,N_1,N_2$ there.

So there are $N_\alpha,f^\alpha_s$ (with $s \in I_\beta$) such that $N_\beta
\le_\gk N_\alpha$ and $s \in I_\beta \Rightarrow f^\alpha_s \circ
g_s = h_s$.  Lastly, for $s \in I_\alpha \setminus I_\beta$ we choose
$h_s = f^\alpha_s \circ g_s$, so we are clearly done.

\sn
2) Similarly, noting that in the last case, the result has cardinality $\leq \chi$ by \ref{f19}(2) or \ref{f28}.
\end{PROOF}

\bigskip

\subsection {Basic Notions} \label{5B}\
\bigskip

As in \cite[\S1]{Sh:600}, we now recall the definition of orbital types
(note that it is natural to look at types only over models which are
amalgamation bases recalling \textbf{Ax.IV}(b) implies every $M \in K_\gk$ is).

\begin{definition}\label{f31}
1) For $\chi \in [\lambda_\gk,\mu_\gk)$ and 
$M \in K_\chi$ we define $\clS(M)$ as 
$$\{ \ortp(a,M,N) : M \le_\gk N \in K_{\leq\chi^{<\kappa}} \text{ and } a \in N\},$$ 
where $\ortp(a,M,N) = (M,N,a)/\cE_M$, where $\cE_M$ is the 
transitive closure of $\cE^\at_M$, and the 
two-place relation $\cE^\at_M$ is defined as follows.

\begin{enumerate}
    \item[$\circledast$] $(M,N_1,a_1)\ \cE^\at_M\  (M,N_2,a_2)$ \underline{iff}:
    \begin{enumerate}
        \item $M \le_\gk N_\ell$ and $a_\ell \in N_\ell$ for $\ell=1,2$.

        \item $\|M\| \le \|N_\ell\| \leq \chi^{<\kappa}$ for $\ell=1,2$.

        \item There exists an $N \in K_{\leq\chi^{<\kappa}}$ and $\le_\gk$-embeddings 
        $f_\ell : N_\ell \to N$ for $\ell = 1,2$ such that 
        $f_1 \rest M = \id_M = f_2 \rest M$ and $f_1(a_1) = f_2(a_2)$.
    \end{enumerate}
\end{enumerate}

\mn
2) We say ``$a$ realizes $p$ in $N$" for $a \in N$ and $p \in \clS(M)$  
\underline{when} (letting $\chi = \|M\|$) for some $N'$ we have 
$M \le_\gk N' \le_\gk N$, $a \in N'$, and $p = \ortp(a,M,N')$. So necessarily
$M,N' \in K_{\chi^{<\kappa}}$, but possibly $N \notin K_{\leq\chi^{<\kappa}}$. 

\sn
3) We say ``$a_2$ strongly\footnote{Note that $\cE^\at_M$ is not an
equivalence relation, and $\cE_M$ certainly isn't, in general.}
realizes $(M,N^1,a_1)/\cE^\at_M$ in $N$" \underline{when} for some
$N^2$ we have $M \le_\gk N^2 \le_\gk N$ and
$a_2 \in N^2$ and $(M,N^1,a_1)\, \cE^\at_M \,(M,N^2,a_2)$. 

\sn
4) We say $M_0$ is a $\le_{\gk[\chi_0,\chi_1)}$-amalgamation
base if this holds in $\gk_{[\chi_0,\chi_1)}$; see below.

\sn
4A) We say $M_0 \in \gk$ is an amalgamation base or 
$\le_\gk$-amalgamation base \underline{when}: for every
$M_1,M_2 \in \gk$ and $\le_\gk$-embeddings
$f_\ell:M_0 \to M_\ell$ (for $\ell = 1,2$) there is $M_3 \in
\gk_\lambda$ and 
$\le_\gk$-embeddings $g_\ell:M_\ell \to
M_3$ (for $\ell=1,2$) such that $g_1 \circ f_1 = g_2 \circ f_2$. 

\sn
5) We say $\gk$ is stable in $\chi$ \underline{when}:
\mn
\begin{enumerate}[(a)]
    \item  $\lambda_\gk \le \chi < \mu_\gk$
\sn
    \item   $M \in K_\chi \Rightarrow |\clS(M)| \le \chi$ 
\sn
    \item  
    $\chi = \chi^{< \kappa}$ 
\sn
    \item $\gk_\chi$ has amalgamation.
\end{enumerate}
\mn
6) We say $p = q \rest M$ if $p \in \clS(M)$, $q \in \clS(N)$, $M \le_\gk N$, and for some $N^+$ such that $N \le_\gk N^+$ and $a \in N^+$ we
have $p = \ortp(a,M,N^+)$ and $q = \ortp(a,N,N^+)$. Note that $p
\rest M$ is well defined if $M \le_\gk N$ and $p \in \clS(N)$. 

\sn
7) For finite $m$, for $M \le_\gk N$ and $\bar a \in {}^m\! N$, we can define
$\ortp(\bar a,N,N)$ and $\clS^m(M)$ similarly, and let $\clS^{< \omega}(M)
= \bigcup\limits_{m < \omega} \clS^m(M)$. (But we shall not use this in any
essential way, hence we choose $\clS(M) = \clS^1(M)$.)
\end{definition}

\begin{remark}\label{f32}
We may replace \ref{f31}(5)(c) by
\begin{enumerate}
    \item[(c)$'$] $\chi \in \Car_\gk$, which means $\chi = \chi^{< \kappa}$ or at least the conclusion of \ref{f19} holds.
\end{enumerate}
In \ref{f34}(1) we change the default value of $\chi$ to 
$\chi = \rnd_\gk(\|N\|)$ (where 
$\rnd_\gk(\theta) \defeq \min(\Car_\gk \setminus \theta)$) so it is 
$\le \|N\|^{<\kappa(\gk)}$ (similarly in \ref{f37}(1)).
\end{remark}

\begin{definition}\label{f34}
1) We say $N$ is $\chi$-\emph{universal above or over} $M$ \underline{when} 
$\chi \in [\lambda_\gk,\mu_\gk)$, $M \in K_{\leq\chi}$, and for every $M'$ with 
$M \le_\gk M' \in K^\gk_\chi$, there is a $\le_\gk$-embedding of 
$M'$ into $N$ over $M$.  If we omit $\chi$ we mean $\|N\|^{<\kappa(\gk)}$; clearly this implies that $M$ is a 
$\le_{\gk_{[\chi_0,\chi_1]}}$-amalgamation base, where $\chi_0 = \|M\|$ 
and $\chi_1 = \|N\|^{<\kappa}$.
 
\sn
2)  $K^3_\gk = \{(M,N,a):M \le_\gk N,\ a \in N \setminus M$ and
$M,N \in K_\gk\}$, with the partial order ${\le} = {\le_\gk}$ defined by
$(M,N,a) \le (M',N',a')$ iff $a = a'$, $M \le_\gk M'$ and $N
\le_\gk N'$.  

\sn
3) We say $(M,N,a)$ is minimal if $(M,N,a) \le (M',N_\ell,a)
\in K^3_\gk$ for $\ell =1,2$ implies $\ortp(a,M',N_1) = 
\ortp(a,M',N_2)$ and moreover, $(M',N_1,a)\ \cE^\at_\lambda\ (M',N_2,a)$
(this is not needed if every $M' \in K_\lambda$ is an amalgamation basis).

\sn
4) $K^{3,\gk}_\lambda$ is defined similarly using
$\gk_{[\lambda,\rnd_\gk(\lambda)]}$.  
\end{definition}

\noindent
Generalizing superlimit, we have more than one reasonable choice.

\begin{definition}\label{f37}
1) For $\ell=1,2$ and $\chi = \chi^{<\kappa} \in [\lambda_\gk,\mu_\gk)$ we say $M^* \in K^\gk_\chi$ is 
{superlimit}$_\ell$ (or $(\chi,{\ge \kappa})$-superlimit$_\ell$)
\underline{when}: (we may omit $\ell$ in the case $\ell=2$) 
\mn
\begin{enumerate}
    \item[(a)]   it is universal, (i.e., every $M \in K^\gk_\chi$ 
    can be properly $\le_\gk$-embedded into $M^*$), and
\sn
    \item[(b)]   \textbf{Case 1}:  $\ell=1$. If  $\LL M_i : i \le \delta \RR$ 
    is $\le_\gk$-increasing, $\cf(\delta) \ge \kappa$, $\delta < \chi^+$, 
    and $i < \delta \Rightarrow M_i \cong M^*$ then $M_\delta \cong M^*$.
\sn
    \item[]   \textbf{Case 2}:  $\ell=2$. If $I$ is a $({<} \kappa)$-directed 
    partial order of cardinality $\le \chi$, $\LL M_t : t \in I\RR$ is 
    $\le_\gk$-increasing, and $t \in I \Rightarrow M_t \cong M^*$ then 
    $\bigcup\{M_t : t \in I\} \cong M^*$.
\end{enumerate}
\mn
2) $M$ is $\chi$-saturated above $\theta$ \underline{if} 
$\|M\| \ge \chi > \theta \ge \LST(\gk)$, and also $N \le_\gk M$, $\theta \le \|N\| < \chi$, 
$p \in \clS_\gk(N)$ imply $p$ is strongly realized in $M$.  Let ``$M$ is
$\chi^+$-saturated" mean that ``$M$ is $\chi^+$-saturated above
$\chi$." Let $$K(\chi^+\text{-saturated}) = \{M \in K : M \text{ is $\chi^+$-saturated}\}$$
 and let ``$M$ is saturated" mean ``$M$ is $\|M\|$-saturated above some 
 $\theta < \|M\|$". 
\end{definition}


\begin{definition}\label{f40}
1) We say $N$ is $(\chi,\sigma)$-\emph{brimmed} over $M$ \underline{when} we
can find a sequence $\LL M_i : i < \sigma \RR$ which is $\le_\gk$-increasing
semi-continuous, $M_i \in K_\chi$, $M_0 = M$, $M_{i+1}$ is 
$\le_\gk$-universal over $M_i$, and $\bigcup\limits_{i < \sigma} M_i = N$.  
We say $N$ is $(\chi,\sigma)$-brimmed over $A$ if $A \subseteq N \in K_\chi$ 
and we can find $\LL M_i : i < \sigma \RR$ as in part (1) such that $A \subseteq M_0$; 
if $A = \varnothing$ we may omit ``over $A$".

\sn
2) We say $N$ is $(\chi,*)$-brimmed over $M$ if for every $\sigma
\in [\kappa,\chi)$, $N$ is $(\chi,\sigma)$-brimmed over $M$.  We say $N$ is
$(\chi,*)$-brimmed if $N$ is $(\chi,*)$-brimmed over $M$ for some $M$. 

\sn
3) If $\alpha < \chi^+$, let ``$N$ is $(\chi,\alpha)$-brimmed over $M$" 
mean $M \le_\gk N$ are from $K_\chi$ and $\cf(\alpha) \ge \kappa \Rightarrow N$ is $(\chi,\cf(\alpha))$-brimmed over $M$.
\end{definition}

\sn
Recall

\begin{claim}\label{f43}
1) If $\gk$ is a DAEC (or just 0-DAEC with amalgamation), stable in 
$\chi$, and $\sigma = \cf(\sigma)$ (so $\chi \in [\lambda_\gk,\mu_\gk)$) 
\underline{then} for every $M \in K^\gk_\chi$ there is 
$N \in K^\gk_\chi$ universal over $M$
which is $(\chi,\sigma)$-brimmed over $M$ 
(hence is $S^\chi_\sigma$-limit: see \cite{Sh:88r}, not used).

\sn
2) If $N_\ell$ is $(\chi,\theta)$-brimmed over $M$ for $\ell =1,2$ and 
$\kappa \le \theta = \cf(\theta) \le \chi^+$
\underline{then} $N_1,N_2$ are isomorphic over $M$. 

\sn
3) If $M_2$ is $(\chi,\theta)$-brimmed over $M$ and $M_0 \le_\gs M_1$ 
then $M_2$ is $(\chi,\theta)$-brimmed over $M_0$.
\end{claim}

\begin{PROOF}{\ref{f43}} 
Straightforward for part (1): recall clause (c) of Definition
\ref{f43}(5).  

\sn
2),3) As in \cite{Sh:600}.
\end{PROOF}

\bn
\centerline {$* \qquad * \qquad *$}
\bigskip

\subsection {Lifting such classes to higher cardinals}\label{5C}

Here we deal with lifting; there are two aspects.  First, if
$\gk^1,\gk^2$ agree in $\lambda$ they agree in every higher cardinal.
Second, given $\gk$ we can find $\gk_1$ with $\mu_{\gk_1} =
\infty$ and $(\gk_1)_\lambda = \gk_\lambda$.

\begin{theorem}\label{f46}
1) If $\gk^\ell$ is a $(\mu,\lambda,\kappa)$-AEC for $\ell=1,2$ 
and $\gk^1_\lambda = \gk^2_\lambda$ then $\gk^1 = \gk^2$.

\sn
2) If $\gk_\ell$ is a $(\mu_\ell,\lambda,\kappa)$-DAEC for $\ell=1,2$ and $\gk^1$ 
satisfies \textbf{Ax.IV}(d), $\mu_1 \le \mu_2$, and
$\gk^1_\lambda = \gk^2_\lambda$ \underline{then} $\gk_1 = \gk_2[\lambda,\mu_1)$.
\end{theorem}

\begin{PROOF}{\ref{f46}}
By \ref{f27}.
\end{PROOF}

\begin{theorem}\label{f49}
\emph{\textbf{The lifting-up Theorem}}

\sn
1) If $\gk_\lambda$ is a $(\lambda^+,\lambda,\kappa)$-DAEC$^\pm$ \underline{then} the pair $\gk' = (K',\le_{\gk'})$ defined below is an 
$(\infty,\lambda,\kappa)$-DAEC$^\pm$, where
\mn
\begin{enumerate}[$(A)$]
    \item   $K'$ is the class of $M$ such that $M$ is a $\tau_{\gk_\lambda}$-model, and for some $I$ and $\olsi M$ we have:
    \begin{enumerate}[$(a)$]
        \item   $I$ is a $\kappa$-directed partial order,
\sn
        \item   $\olsi M = \LL M_s : s \in I \RR$,
\sn
        \item $M_s \in K_\lambda$,
\sn
        \item $I \models ``s < t" \Rightarrow M_s \le_{\gk_\lambda} M_t$.
\sn
        \item  If $J \subseteq I$ has cardinality $\le \lambda$ and is $\kappa$-directed, and $M_J$ is the union of $\LL M_t : t \in J\RR$ (in the sense of \textbf{Ax.III}(d) from Definition \ref{f5}) \underline{then} $M_J$ is a submodel of $M$.
\sn
        \item $M = \bigcup\{M_J:J \subseteq I$ is $\kappa$-directed of cardinality $\le \lambda\}$ (in the sense of \textbf{Ax.IV}(d) of Definition \ref{f5}). 
    \end{enumerate}
\sn
    \item[$(A)'$]  We call such $\LL M_s:s \in I \RR$ a \emph{witness} for $M \in K'$, and we call it \emph{reasonable} if $|I| \le \|M\|^{< \kappa}$.
\sn
    \item[$(B)$]   $M \le_{\gk'} N$ \underline{iff} for some $I,J,\olsi M$ we have:
    \begin{enumerate}
        \item[$(a)$]   $J$ is a $\kappa$-directed partial order,
\sn
        \item[$(b)$]   $I \subseteq J$ is $\kappa$-directed,
\sn
        \item[$(c)$]   $\olsi M =  \LL M_s : s \in J \RR$ and is $\le_{\gk_\lambda}$-increasing,
\sn
        \item[$(d)$]   $\LL M_s:s \in J\RR$ is a witness for $N \in K'$,
\sn
        \item[$(e)$]   $\LL M_s:s \in I\RR$ is a witness for $M \in K'$.
    \end{enumerate}
\sn
    \item[$(B)'$]  We call such $I, \LL M_s:s \in J \RR$ witnesses for $M \le_{\gk'} N$, or say\\ $\big(I,J,\LL M_s:s \in J \RR\big)$ \emph{witnesses}
$M \le_{\gk'} N$. 
\end{enumerate}
\mn
2) If $\gk_\lambda$ satisfies \textbf{Ax.III}(a) then so does $\gk'$. 
\end{theorem}

\begin{PROOF}{\ref{f49}}
1)  Let us check the axioms one by one.

\mn
\textbf{Ax.O(a),(b),(c) and (d):}  $K'$ is a class of 
$\tau_{\gk_\lambda}$-models, $\le_{\gk'}$ a two-place relation on $K$,
$K'$ and $\le_{\gk'}$ are closed under
isomorphisms, and $M \in K' \Rightarrow \|M\| \ge \lambda$, etc.

\noindent
[Why?  trivially.]

\bn
\textbf{Ax.I(a):}  If $M \le_{\gk'} N$ then $M \subseteq N$. 

\noindent
[Why?  We use smoothness for $\kappa$-directed unions; i.e., \textbf{Ax.IV}(x).]



\bn
\textbf{Ax.II(a):}  $M_0 \le_{\gk'} M_1 \le_{\gk'} M_2$ 
implies $M_0 \le_{\gk'} M_2$ and $M \in K' \Rightarrow M \le_{\gk'} M$. 

\noindent
Why?  The second phrase is trivial.  For the first phrase let for $\ell \in
\{1,2\}$ the $\kappa$-directed partial orders $I_\ell \subseteq J_\ell$ and
$\olsi M^\ell = \LL M^\ell_s:s \in J_\ell \RR$ witness 
$M_{\ell-1} \le_{\gk'} M_\ell$.

We first observe
\mn
\begin{enumerate}
    \item[$\odot$] In clause (A)(f) of this theorem, if $J_\bullet \subseteq J_{\bullet\bullet} \subseteq J$ are $({<}\kappa)$-directed \underline{then} $M_{J_\bullet} \leq_{\gk_\lambda} M_{J_{\bullet\bullet}}$.

    \hspace{-1.6cm}[Why? By \textbf{Ax.VI}(d).]
\sn    
    \item[$\boxdot$]   If $I$ is a $\kappa$-directed partial order, 
    $\LL M^\ell_t : t \in I\RR$ is a $\le_{\gk_\lambda}$-increasing 
    sequence witnessing 
    $M_\ell \in K'$ for $\ell=1,2$, and $t \in I \Rightarrow M^1_t \le_{\gk_\lambda} M^2_t$ \underline{then} $M_1 \le_\gk M_2$.
\end{enumerate}
\mn
[Why?  Let $I_1$ be the partial order with set of elements $I \times
\{1\}$ ordered by $(s,1) \le_{I_1} (t,1) \Leftrightarrow s \le_I t$.
Let $I_2$ be the partial order with set of elements $I \times \{1,2\}$
ordered by $(s_1,\ell_1) \le_{I_2} (s_2,\ell_2) \Leftrightarrow s_1
\le_I s_2 \wedge \ell_1 \le \ell_2$.  Clearly $I_1 \subseteq I_2$ are
both $\kappa$-directed.

Let $M_{(s,1)} = M^1_s$ and $M_{(s,2)} = M^2_s$, so clearly $\olsi M =
\LL M_t : t \in I_2\RR$ is a $\le_{\gk_\lambda}$-increasing, $I$-directed
sequence witnessing $M_2 \in K'$, and $(I_1,I_2,\olsi M)$
witnesses $M_1 \le_{\gk'} M_2$, so we have proved $\boxdot$.]

Without loss of generality $J_1,J_2$ are disjoint.  Let $\chi
= \big(|J_1| + |J_2| \big)^{< \kappa}$ so $\lambda \le \chi < \mu_\gk = \infty$, and let

\begin{align*}
\cU \defeq \Big\{u \subseteq J_1 \cup J_2 : &\ |u| \le \lambda,\ u \cap I_\ell \text{ is $\kappa$-directed under $\le_{I_\ell}$ for } \ell = 1,2,  \\
   &\ u \cap J_\ell \text{ is $\kappa$-directed under $\le_{I_\ell}$ for } \ell=1,2, \\
   &\ \text{and } \bigcup\big\{|M^2_t| : t \in u \cap I_2\big\} = \bigcup\big\{|M^1_t| : t \in u \cap J_1\big\}\Big\}.
\end{align*}

\mn
Let $\LL u_\alpha:\alpha < \alpha^*\RR$ list $\cU$, and
we define a partial order $I$ as follows:
\mn
\begin{enumerate}
    \item[(a)$'$]   Its set of elements is $\{\alpha < \alpha^* : 
    \text{for no $\beta < \alpha$ do we have } u_\beta \subseteq u_\alpha\}$.
\sn
    \item[(b)$'$]  For $\alpha,\beta \in I$, $\alpha \le_I \beta$ 
    iff $u_\alpha \subseteq u_\beta$.
\end{enumerate}
\mn
Note that the set $I$ may have cardinality $\big(\sum\limits_{i<\delta}
\|M_i\| \big)^{<\kappa}$ which may be $> \lambda$.

As in the proof of \ref{f19}, $I$ is $\kappa$-directed.

For $\ell = 0,1,2$ and $\alpha \in I$, let $M_{\ell,\alpha}$ be
\mn
\begin{enumerate}
    \item   The $\le_\gk$-union of $\LL M^0_t : t \in u_\alpha \cap I_1\RR$ if $\ell=0$.
\sn
    \item    The $\le_\gk$-union of the $\le_{\gk_\lambda}\!$-directed sequence 
    $\LL M^1_t : t \in J_1\RR$ 
    when $\ell=1$.
\sn
    \item    The $\le_\gk$-union of the $\le_{\gk_\lambda}\!$-directed sequence 
    $\LL M^2_t:t \in J_2\RR$ when $\ell = 2$.
\end{enumerate}
\mn
Now
\mn
\begin{enumerate}
    \item[$(*)_1$]  If $\ell=0,1,2$ and $\alpha \le_I \beta$ then
    $M^\ell_\alpha \le_{\gk_\lambda} M^\ell_\beta$.
\sn
    \item[$(*)_2$]  If $\alpha \in I$ then $M^0_\alpha 
    \le_{\gk_\lambda} M^1_\alpha \le_{\gk_\lambda} M^2_\alpha$.
\sn
    \item[$(*)_3$]  $\LL M_{\ell,\alpha} : \alpha \in I\RR$ is a
    witness for $M_\ell \in K'$.
\sn
    \item[$(*)_4$]  $M_{0,\alpha} \le_{\gk_\lambda}
    M_{2,\alpha}$ for $\alpha \in I$.
\end{enumerate}
\mn
Together by $\boxdot$ we get that $M_0 \le_{\gk'} M_2$ as required.

\bn
\textbf{Ax.III(a):}  In general.

Let $(I_{i,j},J_{i,j},\olsi M^{i,j})$ witness $M_i \le_{\gk'} M_j$ 
when $i \le j < \delta$, and without loss of generality 
$\LL J_{i,j} : i < j < \delta\RR$ are pairwise disjoint.  
Let $\cU$ be the family of sets $u$ such that for some 
$v \in [\delta]^{\le \lambda}$,
\mn
\begin{enumerate}
    \item $v \subseteq \delta$ has cardinality $\le \lambda$ and has order type of cofinality $\ge \kappa$.
\sn
    \item $u \subseteq \bigcup\{J_{i,j} : i < j$ are from $v\}$ has cardinality $\le \lambda$.
\sn
    \item   For $i \leq j$ from $v$, the set $u \cap J_{i,j}$ is $\kappa$-directed under $\le_{J_{i,j}}$ and $u \cap I_{i,j}$ is $\kappa$-directed under $\le_{I_{i,j}}$.
\sn

    \item  If $i(0) \le i(1) \le i(2)$ are from $v$ then
    $$\bigcup\big\{M^{i(0),i(1)}_{s} : s \in u \cap J_{i(0),i(1)}\big\} =
    \bigcup\big\{M^{i(1),i(2)}_s : s \in u \cap I_{i(1),i(2)}\big\}.$$

    \item If $i(0) \le j(1)$ and $i(1) \le j(1)$ are from $v$ then 
    $$\bigcup\big\{M^{i(0),j(1)}_s : s \in u \cap J_{i(0),j(1)}\big\} =
    \bigcup\big\{M^{i(1),j(1)}_s : s \in u \cap J_{i(1),j(1)}\big\}.$$
\end{enumerate}
\mn
Let the rest of the proof be as in the proof of \textbf{Ax.II}(a).

\bn
\textbf{\textbf{Ax.IV}(a):}

Similar, but $\cU = \{u \subseteq I:u$ has cardinality $\le \lambda$
and is $\kappa$-directed$\}$.

\bn
\textbf{\textbf{Ax.III}(d):}

Recall that we are assuming $\gk$ satisfies \textbf{Ax.III}(d).  Similar proof.

\bn
\textbf{\textbf{Ax.IV}(d):}

Again, we are assuming $\gk$ satisfies \textbf{Ax.IV}(d).

\bn
\textbf{Ax.V:}  Assume $N_0 \le_{\gk'} M$ and $N_1 \le_{\gk'} M$.

If $N_0 \subseteq N_1$, then $N_0 \le_{\gk'} N_1$. 

\sn
[Why?  Let $\big(I_\ell,J_\ell,\LL M^\ell_s:s \in J_\ell \RR\big)$ witness 
$N_\ell \le_\gk M$ for $\ell=0,1$; without loss of generality $J_0,J_1$ are disjoint.

Let 
\begin{align*}
\cU \defeq \Big\{ u \subseteq J_0 \cup J_1 : &\ |u| \le \lambda,\ 
    u \cap J_\ell \text{ and } u \cap I_\ell \text{ are } \kappa\text{-directed for } \ell=0,1,\\
  &\text{ and } 
  \bigcup \big\{|M^0_s| : s \in u \cap J_0\} = \bigcup\big\{|M^0_s| : s \in u \in
  J_1 \big\}\Big\}.
\end{align*}

For $u \in \cU$ let
\sn
\begin{itemize}
    \item  $M_u = M \rest \bigcup\{M^\ell_s : s \in u \cap J_\ell\}$ for $i=0,1$.
\sn
    \item $N_{\ell,u} = N_\ell \rest \{M^\ell_s : s \in u \cap I_\ell\}$.
\end{itemize}
\sn
Let
\sn
\begin{enumerate}
    \item[$(*)$] 
    \begin{enumerate}
        \item $(\cU,\subseteq)$ is $\kappa$-directed.

        \item $N_{\ell,u} \le_\gk M$

        \item $M_{\ell,u} \le_\gk M_{\ell,v}$ when $u \subseteq v$ are from $\cU$ and $\ell=0,1$.

        \item $M_{0,u} \le_\gk M_{1,u}$

        \item $N_\ell = \bigcup\{N_{\ell,u} : u \in \cU\}$
    \end{enumerate}
\end{enumerate}
\mn
By $\boxdot$ above we are done.

\bn
\textbf{Ax.VI}:  $\LST(\gk') = \lambda$. 

\sn
[Why?  Let $M \in K'$, $A \subseteq M$, $|A| + \lambda \le \chi < \|M\|$, 
and let $\LL M_s : s \in I \RR$ witness $M \in K'$. 
Without loss of generality $|A| = \chi^{< \kappa}$.  
Now choose a directed $I \subseteq J$ of cardinality 
$\le |A| = \chi^{< \kappa}$ such that 
$A \subseteq M' \defeq \bigcup\limits_{s \in I} M_s$ and so 
$(I,J,\LL M_s : s \in J\RR)$ witnesses $M' \le_{\gk'} M$. So 
as $A \subseteq M'$ and $\|M'\| \le |A| + \mu$ we are done.] 

\bn
2) \textbf{Ax.III}(b):  Assume that $\LL M_i : i < \delta \RR$ 
is a $\le_\gk$-increasing sequence and 
$$\big\| \textstyle\bigcup\{M_i : i < \delta\} \big\| < \mu.$$ 
We have to prove that there is $M \in K$ such that $i < \delta \Rightarrow M_i \leq_{\gk'} M$ (as always, assuming that $\gk_\lambda$ satisfies \textbf{Ax.III}(b)).

If $\cf(\delta) \geq \kappa$, we use \textbf{Ax.III}(a) (which was proved above), so we may assume $\cf(\delta) < \kappa$. By renaming, without loss of generality $\delta < \kappa$, and we continue in the proof of \textbf{Ax.II}(a).
\end{PROOF}

\sn
Also, if two such DAECs have some cardinal in common then we can
put them together.

\begin{claim}\label{f52}
Let $\iota \in \{0,1,2,4\}$, assume $\lambda_1 < \lambda_2 < \lambda_3$, and
\mn
\begin{enumerate}[$(a)$]
    \item   $\gk^1$ is an $(\lambda^+_2,\lambda_1,\kappa)$-$2$-DAEC and 
    $K^1 = K_{\gk^1}$.
\sn
    \item  $\gk^2$ is a $(\lambda_3,\lambda_2,\kappa)$-$\iota$-DAEC.
\sn
    \item  $K^{\gk^1}_{\lambda_2} = K^{\gk^2}_{\lambda_2}$ and 
    ${\le_{\gk^2} \rest K^{\gk^2}_{\lambda_2}} = {\le_{\gk^1} \rest K^{\gk^1}_{\lambda_2}}$.
\sn
    \item  We define $\gk$ as follows: $K_\gk = K_{\gk^1} \cup K_{\gk^2}$, $M \le_\gk N$
    iff $M \le_{\gk^1} N$ or $M \le_{\gk^2} N$ or for some $M'$, 
    $M \le_{\gk^1} M' \le_{\gk^2} N$.
\end{enumerate}
\sn
\underline{Then} $\gk$ is an $(\lambda_3,\lambda_1,\kappa)$-$\iota$-DAEC.
\end{claim}

\begin{PROOF}{\ref{f52}}
Straightforward.  E.g.: 

\mn
\textbf{\textbf{Ax.III}(d):}  $\LL M_s:s \in I\RR$ is a
$\le_\gs$-$\kappa$-directed system.

If $\|M_s\| \ge \lambda_2$ for some $s$, use $\LL M_t : s \le
t \in I\RR$ and clause (b) of the assumption.  If $\bigcup\{M_s : s \in I\}$ 
has cardinality $\le \lambda_2$ use clause (a) in the
assumption.  If neither one of them holds, recall $\lambda_2 =
\lambda^{< \kappa}_2$ by clause (b) of the assumption, and let
\[
\cU = \big\{u \subseteq I : |u| \le \lambda_2,\ u \text{ is }
\kappa\text{-directed (in I), and } \textstyle\bigcup\{M_s : s \in u\} \text{ has
  cardinality } \lambda \big\}.
\]

Easily, $(\cU,\subseteq)$ is $\lambda_2$-directed. For $u \in J$, let
$M_u$ be the $\le_\gs$-union of $\LL M_s:s \in u\RR$.  Now
by clause (a) of the assumption
\mn
\begin{enumerate}
    \item[$(*)_1$]  $M_u \in K^{\gk^1}_{\lambda_2} = K^{\gk^2}_{\lambda_2}$
\sn
    \item[$(*)_2$]  If $u_1 \subseteq v$ are from $\cU$ then 
    $M_u \le_{\gk^1} M_v$, $M_u \le_{\gk^2} M_v$.
\end{enumerate}
\mn
Now use clause (b) of the assumption.

\bn
\textbf{Axiom V:}  We shall freely use
\sn
\begin{enumerate}
    \item[$(*)$]   $\gk^2_{\lambda_2} =\gk^1_{\lambda_2} = \gk_{\lambda_2}$
\end{enumerate}
\sn
So assume $N_0 \le_\gk M$, $N_1 \le_\gk M$, $N_0 \subseteq N_1$. 

Now if $\|N_0\| \ge \lambda_2$ use assumption (b), so we can assume
$\|N_0\| < \lambda_2$. If $\|M\| \le \lambda_2$ we can use assumption (a), 
so assume $\|M\| > \lambda_2$; by the definition of $\le_\gk$ there
is $M'_0 \in K^{\gk^1}_{\lambda_2} = K^{\gk^2}_{\lambda_2}$ such that 
$N_0 \le_{\gk^1} M'_0 \le_{\gk^2} M$.  First assume 
$\|N_1\| \le \lambda_2$, so we can find 
$M'_1 \in K^{\gk^1}_{\lambda_2}$ such that
$N_1 \le_{\gk^1} M'_1 \le_{\gk^2} M$.

\sn
[Why?  If $N_1 \in K^{\gk^1}_{< \lambda_2}$ by the definition of 
$\le_\gk$, and if $N_1 \in K^{\gk^1}_{\lambda_2}$ just 
choose $M'_1 = N_1$.]  

Now we can, by assumption (b), 
find $M'' \in K^{\gk^1}_{\lambda_2}$ such that $M'_0 \cup M'_1 \subseteq M'' 
\le_\gk M$, hence by assumption (b) (i.e. \textbf{Ax.V} for $\gk^2$)
we have $M'_0 \le_\gk M''$, $M'_1 \le_\gk M''$.
As $N_0 \le_\gk M'_0 \le_\gk M'' \in K^\gk_{\le \lambda_2}$ 
by assumption (a) we have $N_0 \le_\gk M''$, and similarly 
we have $N_1 \le_\gk M''$.  So
$N_0 \subseteq N_1$, $N_0 \le_\gk M''$, $N_1 \le_\gk M'$, so by 
assumption (b) we have $N_0 \le_\gk N_1$.

We are left with the case $\|N_1\| > \lambda$. By assumption (b) there is
$N'_1 \in K_{\lambda_2}$ such that $N_0 \subseteq N'_1 \le_{\gk^2} N_2$. 
{Also} by assumption (b), we have $N'_1 \le_\gk M$, so by the previous paragraph
we get $N_0 \le_\gk N'_1$; together with the previous sentence we have
$N_0 \le_{\gk^1} N'_1 \le_{\gk^2} N_1$ so by the definition of
$\le_\gk$ we are done. 
\end{PROOF}

\begin{definition}\label{f55}
If $M \in K_\chi$ is $(\chi,\ge \kappa)$-superlimit$_1$ let
$$K^{[M]}_\chi = \{N \in K_\chi:N \cong M\}$$ 
and $\gk^{[M]}_\chi = (K^{[M]}_\chi,\le_\gk \rest K^{[M]}_\chi)$, 
and $\gk^{[M]}$ is the $\gk'$ we get in \ref{f49}(1), with $(\gk^{[M]},\gk)$ here standing in for $(\gk_\lambda,\gk')$ there. 
\end{definition}

\begin{claim}\label{f58}
1) If $\gk$ is an $(\mu,\lambda,\kappa)$-AEC, $\lambda \le \chi < \mu$, 
$M \in K_\chi$ is $(\chi,{\ge} \kappa)$-superlimit$_1$ (see Definition \ref{f37}) \underline{then}
$\gk^{[M]}_\chi$ is a $(\chi^+,\chi,\kappa)$-DAEC.

\sn
2) If in addition $\gk$ is a $(\mu,\lambda,\kappa)$-DAEC$^\pm$ \underline{then} 
$\gk^{[M]}_\chi$ is a $(\chi^+,\chi,\kappa)$-DAEC$^\pm$.

\sn
3) $[\gk$ satisfies \textbf{Ax.IV}(d).$]$

$M$ is $(\chi,\geq\kappa)$-superlimit$_1$ \underline{iff} $M$ is $(\chi,\geq\kappa)$-superlimit$_2$.
\end{claim}

\begin{PROOF}{\ref{f58}}
Easy.
\end{PROOF}

\newpage

\section {pr frames} \label{6}

Below, the main case is $\iota = 4$.

\begin{definition}\label{g2}
Here $\iota = 0,1,2,3,4$.  We say that $\gs$ is a good
$(\mu,\lambda,\kappa)$-$\iota$-frame \underline{when} $\gs$ 
consists of the following objects
satisfying the following condition: $\mu,\lambda,\kappa$ (so we should
write $\mu_\gs,\lambda_\gs,\kappa_\gs$ but we
may ignore them when defining $\gs$) and
\mn
\begin{enumerate}
    \item $\gk = \gk_\gs$ is a $(\mu,\lambda,\kappa)$-4-DAEC (see \ref{f5}(4)), 
    so we may write $\gs$ instead of $\gk$, e.g. $\le_\gs$-increasing, 
    etc., and $\chi \in [\lambda,\mu) \Rightarrow \LST(\chi^{< \kappa})$.
\sn
    \item   $\gk$ has a $(\lambda,\ge \kappa)$-superlimit model $M^*$ 
    which\footnote{Follows by (C), in fact.} 
    is not $<_\gk$-maximal --- i.e.:
    \begin{enumerate}
        \item  $M^* \in K^\gs_\lambda$
\sn
        \item  If $M_1 \in K^\gs_\lambda$ \underline{then} for some $M_2$, 
        $M_1 <_\gs M_2 \in K^\gs_\lambda$ and $M_2$ is isomorphic to $M^*$.
\sn
        \item If $\LL M_i : i < \delta \RR$ is  $\le_\gs$-increasing, 
        $i < \delta \Rightarrow M_i \cong M$, and $\cf(\delta) \ge \kappa$, $\delta < \lambda^+$ \underline{then} $\bigcup\{M_i : i < \delta\}$ is isomorphic to $M^*$.
    \end{enumerate}

    \item $\gk$ has the amalgamation property, the JEP 
    (joint embedding property), and has no $\le_\gk$-maximal member. 
    If $\iota \ge 1$ then $\gk$ has primes over chains (i.e. \textbf{Ax.III}(f)), 
    and if $\iota \ge 4$, $\gk$ has primes over $\leq_\gs$-directed sequences (i.e. \textbf{Ax.IV}(f)). 
\sn
    \item  
    \begin{enumerate}
        \item $\clS^\bs = \clS^\bs_\gs$ 
        (the class of basic types for $\gk_\gs$) is included in \\
        $\bigcup\{\clS(M) : M \in K_\gs\}$ and is closed under isomorphisms including automorphisms. For $M \in K_\lambda$, let 
        $\clS^\bs_\gs(M) = \clS^\bs_\gs \cap \clS(M)$; no harm in allowing types of finite sequences.
\sn
        \item If $p \in \clS^\bs_\gs(M)$, \underline{then} $p$ is non-algebraic (i.e., not realized by any $a \in M$).
        
        \item \textbf{Density}: 
        
        If $M \le_\gk N$ are from $K_\gs$ and $M \ne N$, \underline{then} for some $a \in N \setminus M$  we have $\ortp(a,M,N) \in \clS^\bs$. The intention is that examples are minimal types in \cite{Sh:576} (i.e. \cite{Sh:E46}) and regular types for superstable theories.
\sn
        \item \textbf{bs-Stability}: 
        $\clS^\bs(M)$ has cardinality $\le \|M\|^{< \kappa}$ for $M \in K_\gs$. 
    \end{enumerate}
\sn
    \item
    \begin{enumerate}
        \item $\nonfork{}{}_{} = \nonfork{}{}_\gs$ is a four place relation  called \emph{non-forking}, with $\nonfork{}{}_{}(M_0,M_1,a,M_3)$ implying $M_0 \le_\gk M_1 \le_\gk M_3$ are from $K_\gs$, 
        $a \in M_3 \setminus M_1$, $\ortp(a,M_0,M_3) \in \clS^\bs_\gs(M_0)$, and 
        $\ortp(a,M_1,M_3) \in \clS^\bs(M_1)$. Also, $\nonfork{}{}_{}$ is preserved under isomorphisms. 
        
        We may also write $\nonforkin{M_1}{a}_{M_0}^{M_3}$ and demand that if $M_0 = M_1 \le_\gk M_3$ are both in $K_\lambda$ \underline{then} $\nonfork{}{}_{}(M_0,M_1,a,M_3)$ is equivalent to 
        ``$\ortp(a,M_0,M_3) \in \clS^\bs(M_0)$".  We may state $\nonforkin{M_1}{a}_{M_0}^{M_3}$ as ``$\ortp(a,M_1,M_3)$ does not fork over $M_0$ (inside $M_3$)." 
        (This is justified by clause (b) below.)
        
        [Explanation: The intention is to axiomatize non-forking of types, 
        but we allow ourselves to deal only with basic types. 
        Note that in \cite{Sh:576} (i.e. \cite{Sh:E46}) we know 
        something on minimal types but other types are something else.]

        \item \textbf{Monotonicity:} 
        $$\text{If } M_0 \le_\gk M'_0 \le_\gk M'_1 \le_\gk M_1 \le_\gk M_3 \le_\gk M'_3$$ 
        and $M_1 \cup \{a\} \subseteq M''_3 \le_\gk M'_3$, with all of them in $K_\lambda$, \underline{then}
        $$\nonfork{}{}_{}(M_0,M_1,a,M_3) \Rightarrow \nonfork{}{}_{}(M'_0,M'_1,a,M'_3) \Leftrightarrow \nonfork{}{}_{}(M'_0,M'_1,a,M''_3)$$ 
        \underline{so} it is legitimate to just say ``$\ortp(a,M_1,M_3)$ does not fork over $M_0$". 
        
\sn
        [Explanation: non-forking is preserved by decreasing the type, increasing the basis (i.e. the set over which it does not fork) and increasing or decreasing the model inside which all this occurs.  The same holds for stable theories, only here we restrict ourselves to ``legitimate" types.]
\sn
        \item \textbf{Local Character:} 

        \noindent\textbf{Case 1}: $\iota=1,2,3$.

        If $\LL M_i : i \le \delta\RR$ is $\le_\gs$-semi-continuous, $p \in \clS^\bs(M_\delta)$, and $\cf(\delta) \ge \kappa$ \underline{then} for every $\alpha < \delta$ large enough, $p$ does not fork over $M_\alpha$.

        \sn
        \textbf{Case 2}:  $\iota=4$.  

        If $I$ is a $\kappa$-directed partial order, $\olsi M = \LL M_t : t \in I\RR$ is a 
        $\le_\gs$-directed system, $M$ is its $\le_\gk$-union, $M \le_\gs N$, and
        $\ortp(a,M,N) \in \clS^\bs(M_\delta)$ \underline{then} for every $s \in I$ large enough $\ortp(a,M,N)$ does not fork over $M_s$. 

        \sn
        \textbf{Case 3}:  $\iota=0$.

        Like Case 1, using $(\geq\!\kappa)$-continuity.

        \sn 
        [Explanation: This is a replacement for $\kappa \ge \kappa_r(T)$: if $p \in \clS(A)$ then there is a $B \subseteq A$ of cardinality $< \kappa$ such that $p$ does not fork over $A$.]  
\sn
        \item \textbf{Transitivity}: 

        If $M_0 \le_\gk M'_0 \le_\gk M''_0 \le_\gk M_3$ and $a \in M_3$ and $\ortp(a,M''_0,M_3)$ does not fork over $M'_0$ and $\ortp(a,M'_0,M_3)$ does not fork over $M_0$ (all models are in $K_\lambda$, of course, and necessarily the three relevant types are in $\clS^\bs$), \underline{then} $\ortp(a,M''_0,M_3)$ does not fork over $M_0$.
\sn
        \item \textbf{Uniqueness}: 

        If $p,q \in \clS^\bs(M_1)$ do not fork over $M_0 \le_\gk M_1$ (all in $K_\gs$) and 
        $p \rest M_0 = q \rest M_0$ \underline{then} $p = q$.
\sn
        \item \textbf{symmetry}:  

        \sn
        \textbf{Case 1}: $\iota \ge 3$.

        If $M_0 \le_\gs M_\ell \le_\gs M_3$ and $(M_0,M_\ell,a_\ell) \in K^{3,\pr}_\gs$ (see clause (j) below) for $\ell=1,2$ \underline{then} $\ortp_\gs(a_2,M_1,M_3)$ does not fork over $M_0$ \underline{iff} $\ortp_\gs(a_1,M_2,M_3)$ does not fork over $M_0$.

        \sn
        \textbf{Case 2}:  $\iota=0,1,2$.

        If $M_0 \le_\gk M_3$ are in $\gk_\lambda$, and for $\ell = 1,2$ we have $a_\ell \in M_3$ and $\ortp(a_\ell,M_0,M_3) \in \clS^\bs(M_0)$, \underline{then} the following are equivalent:
        \begin{enumerate}
            \item[$(\alpha)$]  There are $M_1,M'_3$ in $K_\gs$ such that
            $M_0 \le_\gk M_1 \le_\gK M'_3$,  $a_1 \in M_1$, $M_3 \le_\gk M'_3$ and $\ortp(a_2,M_1,M'_3)$ does not fork over $M_0$.
\sn
            \item[$(\beta)$] There are $M_2,M'_3$ in $K_\lambda$ such that
            $M_0 \le_\gk M_2 \le_\gk M'_3$, $a_2 \in M_2$, $M_3 \le_\gk M'_3$ and $\ortp(a_1,M_2,M'_3)$ does not fork over $M_0$. 

            [Explanation: this is a replacement to ``$\ortp(a_1,M_0 \cup
            \{a_2\},M_3)$ forks over $M_0$ iff $\ortp(a_2,M_0 \cup \{a_1\},M_3)$
            forks over $M_0$," which is not well defined in our context.]
        \end{enumerate}
\sn
        \item \textbf{Existence}: 

        If $M \le_\gs N$ and $p \in \clS^\bs(M)$ \underline{then} there is $q \in \clS^\bs(N)$ which is a non-forking extension of $p$.
\sn
        \item \textbf{Continuity}: 
        
        \hspace{-.5cm}\textbf{Case 1:}  $\iota = 1,2,3$.

        If $\LL M_\alpha : \alpha \le \delta\RR$ is $\le_\gs$-increasing and $\le_\gs$-semi-continuous, $M_\delta = \bigcup\limits_{\alpha < \delta} M_\alpha$ (which holds if $\cf(\delta) \ge \kappa$), $p \in \clS_\gs(M_\delta)$, and $p \rest M_\alpha$ does not fork over $M_0$ for $\alpha < \delta$ \underline{then} $p \in \clS^\bs_\gs(M_\delta)$ and it does not fork over $M_0$.

        \mn
        \hspace{-.5cm}\textbf{Case 2:}  $\iota=4$.

        Similarly, but for $\olsi M = \LL M_t : t \in I\RR$, $I$ directed, 
        and\\ $M = \bigcup\{M_t : t \in I\}$ is a $\le_\gs$-upper bound of 
        $\olsi M$.

        \mn
        \hspace{-.5cm}\textbf{Case 3:}  $\iota=0$.

        Like Case 1 for $\olsi M$ $(\geq\!\kappa)$-continuous.
        \begin{enumerate}
            \item If $\iota \geq 1$, $\gs$ has $K^{3,\pr}_\gs$-primes (see \ref{g26} below).
\sn
            \item If $p \in \clS^\bs_\gs(N)$ \underline{then} $p$ does not fork over $M$ for some $M \le_\gs N$ from $K_\lambda$.    
        \end{enumerate}  
\sn
        \item \textbf{Strong continuity}:
        
        \hspace{-.5cm}\textbf{Case 1:}  $\iota = 1,2,3$.

        We have that $\ortp(b,M_\delta,M_{\delta+1})$ does not fork over $M_0$, where
        \begin{enumerate}
            \item $\olsi M = \LL M_i : i \leq \delta + 1\RR$ is $\leq_\gs$-increasing.

            \item $M_\delta$ is prime over $\olsi M \rest \delta$.

            \item $b \in M_{\delta+1} \setminus M_\delta$

            \item $\ortp(b,M_i,M_{\delta+1})$ does not fork over $M_0$ for 
            $i < \delta$.
        \end{enumerate}       

        \mn
        \hspace{-.5cm}\textbf{Case 2:} $\iota = 4,5$.
        
        We have that $\ortp(b,M_\delta,M_{\delta+1})$ does not fork over $M_0$, where
        \begin{enumerate}
            \item $\olsi M = \LL M_s : s \in I\RR$ is $\leq_\gs$-increasing, $I$ a partial order with $0 \in I$ minimal.

            \item $N_0$ is prime over $\olsi M$.

            \item $b \in N_1 \setminus N_0$, where $N_0 \leq_\gs N_1$.

            \item $\ortp(b,M_s,N_1)$ does not fork over $M_0$ for all $s \in I$.
        \end{enumerate}      
    \end{enumerate}  
\end{enumerate}
\end{definition}


\begin{claim}\label{g8}
1) If $\LL M_i : i < \delta \RR$ is $\le_\gk$-increasing, 
$\big(\sum\big\{\|M_i\|:i < \delta\big\}\big) < \mu$, 
$p_i \in \clS^\bs_\gs(M_i)$ does not fork over $M_0$ for $i < \delta$, 
and $i < j \Rightarrow p_j \rest M_i = p_i$ \underline{then}:
\mn
\begin{enumerate}[$(a)$]
    \item  We can find $M_\delta$ such that 
    $i < \delta \Rightarrow M_i \le_\gk M_\delta$.
\sn
    \item   For any such $M_\delta$, we can find $p_\delta \in \clS_\gs(M_\delta)$ 
    such that $\bigwedge\limits_{i < \delta} {p_\delta} \rest M_i = p_i$ 
    and ${p_\delta}$ does not fork over $M_0$.
\sn
    \item   In clause (b), {$p_\delta$} is unique. 
\sn
    \item   If $\ell \ge \kappa \wedge \cf(\delta) \ge \kappa$, we can add 
    $M = \bigcup\{M_\alpha : \alpha < \delta\}$.
\end{enumerate}
\mn
2) Similarly for $\olsi M = \LL M_t : t \in I\RR$, $I$ directed.
\end{claim}

\begin{PROOF}{\ref{g8}}
1) First, choose $M_\delta$ by \ref{g2}, Clause (A).
Second, choose $p_\delta \in \clS^\bs_\gs(M_\delta)$, a
non-forking extension of $p_0$, which exists by Axiom (g) of \ref{g2}(E).   
Now $p_\delta \rest M_i \in \clS^\bs_\gs(M_i)$ does 
not fork over $M_0$ by \ref{g2}(E)(b) and it extends $p_0$, so it is equal to $p_i$ by (E)(e).
Third, $p_\delta$ is unique by (E)(e).

\noindent
2) Should be clear, too.  
\end{PROOF}

\begin{definition}\label{g11}
1) Assume $M_\ell \le_\gs N$ and $p_\ell \in \clS^\bs_\gs (M_\ell)$ 
for $\ell = 1,2$.  We say that
$p_1,p_2$ are parallel \underline{when} some 
$p \in \clS^\bs_\gs(N)$ is a non-forking 
extension of $p_\ell$ for $\ell=1,2$.

\noindent
2) We say $\gs$ is \emph{type-full} when $\clS^\bs_\gs
(M) = \clS^\na_{\gk_\gs}(M)$ for $M \in K_\gs$.

\noindent
3) We say $p \in \clS^\bs_\gs(M)$ is based on $\bar\bfa$ when:
\begin{enumerate}
    \item  $\bar\bfa$ is a sequence from $M$.
\sn
    \item  If $M \le_\gs N$, $q \in \clS^\bs_\gs(N)$ is a non-forking extension of $p$, and $\pi$ is an automorphism of $N$ over $\bar\bfa$ then 
    $\pi(q) = q$. (See \cite{Sh:839} for how we can guarantee there is such $\bar\bfa \in {}^\lambda\!M$, and even $\bar\bfa \in {}^1\!M$.)
\end{enumerate}
\mn
3A) Similarly for $p \in \clS_\gs^\eps(M)$; similarly for part (4).

\noindent
4) We say $\gs$ is $(< \theta)$-based when in clause \ref{g11}(3) above there is such 
$\bar\bfa\in {}^{\theta >}\!M$.
\end{definition}

\begin{definition}\label{g14}
1) We say that NF is a non-forking relation on a
$(\mu,\lambda,\kappa)$-1-DAEC $\gk$ \underline{when}, 
in addition to \ref{g2}(A)-(C):
\mn
\begin{enumerate}
    \item[(F)]  
    \begin{enumerate}
        \item NF is a four-place relation on $\gk_\gs$, and $\NF_\gs(M_0,M_1,M_2,M_3)$ implies $M_0 \le_\gk M_\ell \le_\gk M_1$ and $\NF_\gs$ is preserved by isomorphisms.
\sn
        \item[(b)$_1$] \textbf{Monotonicity}: 

        If $\NF_\gs(M_0,M_1,M_2,M_3)$, $M_0\le_\gs M'_\ell \le_\gs M_\ell$ for $\ell=1,2$, and\\ $M'_1 \cup M'_2 \subseteq M'_3 \le_\gs M$ 
        then $\NF_\gs(M_0,M'_1,M'_2,M'_3)$.
\sn
        \item[(c)] \textbf{Symmetry}: 
        
        $\NF_\gs(M_0,M_1,M_2,M_3)$ implies $\NF_\gs(M_0,M_2,M_1,M_3)$.
\sn
        \item[(d)$_1$] \textbf{Transitivity}: 
        
        If $\NF_\gs(M_{2 \ell},M_{2\ell+1},M_{2\ell + 3},M_{2\ell +4})$ 
        for $\ell = 0,1$ then $\NF_\gs(M_0,M_1,M_4,M_5)$.
        \item[(d)$_2$]  \textbf{Long transitivity}: 
        
        If $\LL (N_i,M_i) : i < \delta \RR$ is an $\NF_\gs$-sequence (i.e.  $M_i$ is $\le_\gs$-increasing, $N_i$ is $\le_\gs$-increasing, $M_i \le_\gs N_i$, $i < j < \delta \Rightarrow \NF_\gs(M_i,N_i,M_j,N_j)$, and 
        $\sum\{\|N_i\|:i < \delta\} < \mu$) \underline{then} we can find $(N_\delta,M_\delta)$ such that $\LL (M_i,N_i) : i \le \delta \RR$ is an $\NF$-sequence. 
\sn
        \item[(d)$_2^{\!+}$] Like (d)$_2$, for directed systems. 
\sn
        \item[(d)$_3$] Moreover, in (d)$_2$, if $\LL M_i : i < \delta\RR$ and $\LL N_i : i < \delta\RR$ are $(\geq\!\kappa)$-continuous and $\cf(\delta) \geq \kappa_\gs$ \underline{then} we can demand $M_\delta = \bigcup\{ M_i : i < \delta\}$, $N_\delta = \bigcup\{ N_i : i < \delta\}$.
    \end{enumerate}
\end{enumerate}
\end{definition}

\begin{definition}\label{g17}
1) Let $\gs$ be a good $\lambda$-frame and
NF a non-forking relation on $\gk$.  We say NF respects $\gs$
\underline{when}: if $\NF_\gs(M_0,M_1,M_2,M_3)$ and $a \in M_2$, 
$\ortp_\gs(a,M_0,M_3) \in \clS^\bs_\gs(M_0)$ \underline{then}
$\ortp_\gs(a,M_1,M_3)$ is a non-forking extension of 
$\ortp_\gs(a,M_0,M_2)$.

\noindent
2) We say $\gs$ is a good $(\lambda,\mu,\kappa)$-$\NF$-frame \underline{when}
it is a good $(\lambda,\mu,\chi)$-frame and $\NF_\gs$ is a
non-forking relation on $\gk_\gs$ which respects $\gs$.
\end{definition}

\begin{definition}\label{g23}
We say that $\gs$ is a \emph{very good} $(\mu,\lambda,\kappa)$-$\NF$-frame
\underline{when} it is a good $(\mu,\lambda,\kappa)$-$\NF$-frame and
\mn
\begin{enumerate}
    \item[(G)]  
    \begin{enumerate}
        \item If $\NF_\gs(M_0,M_1,M_2,M_3)$ \underline{then} there is 
        $M^*_3 \le_\gs M_3$ which is prime over $M_1 \cup M_2$. That is,
        \begin{itemize}
            \item  If $\NF_\gs(M'_0,M'_1,M'_2,M'_3)$ and $f_\ell$ is an isomorphism from $M_\ell$ onto $M'_\ell$ for $\ell =0,1,2$ such that $f_0 \subseteq f_1$ and $f_0 \subseteq f_2$ \underline{then} there is a $\le_\gs$-embedding $f_3$ of $M^*_3$ into $M'_3$ extending 
            $f_1 \cup f_2$.
        \end{itemize}
\sn
        \item $\gk_\gs$ has $K_\gs^{3,\pr}$-primes (see \ref{g26}(3) below).
    \end{enumerate}
\end{enumerate}
\end{definition}

\begin{definition}\label{g26}
0) $K^{3,\bs}_\gs = \{(M,N,a) : M \le_\gs N,\ a \in N$, and 
$\ortp_\gs(a,M,N) \in \clS^\bs_\gs(M)\}$.

\noindent
1) $K^{3,\pr}_\gs = \{(M,N,a) \in K^{3,\bs}_\gs: 
\text{if } M \le N',\ a' \in N',\ \ortp_\gs(a',M,N') 
= \ortp(a,M,N)$ \underline{then} there is a $\le_\gk$-embedding of
$N$ into $N'$ extending id$_M$ and mapping $a$ to $a'\}$.

\noindent
2) $\gk_\gs$ has $K^{3,\pr}_\gs$-primes \underline{if}, for 
every $M \in K_\gs$ and $p \in \clS^\bs_\gs(M)$, 
there are $(N,a)$ such that $(M,N,a) \in K^{3,\pr}_\gs$ 
and $\ortp_\gs(a,M,N) = p$. 
\end{definition}

\begin{definition}\label{g29} [$\iota \ge 3$]

\noindent 
1) Assume $p_1,p_2 \in \clS^\bs(M)$.  We say $p_1,p_2$ are
weakly orthogonal (and denote it $p_1 {\underset \wk \perp} p_2$) when 
the following implication holds: if $M_0 \le_\gs M_\ell \le_\gs M_3$, $(M_0,M_\ell,a_\ell) \in K^{3,\pr}_\gs$, and 
$\ortp_\gs(a_\ell,M_0,M_\ell) = p_\ell$ for $\ell=1,2$ \underline{then} $\ortp_\gs(a_2,M_1,M_3)$ does not
fork over $M_0$ (this is symmetric by Axiom (f) of \ref{g2}(E)).

\noindent
2) We say $p_1,p_2$ are orthogonal (denoted $p_1 \perp p_2$) when: if
$M \le_\gs M_2$, $M_1 \le_\gs M_2$ and $q_\ell \in \clS^\bs(M_2)$
is a non-forking extension of $p_\ell$ and
$q_\ell$ does not fork over $M_1$ \underline{then} $q_1 {\underset \wk \perp} q_2$.

\noindent
3) We say that $\{a_t : t \in I\}$ is independent in
$(M_0,M_1,M_2)$ when:
\mn
\begin{enumerate}
    \item  $a_t \in M_2 \setminus M_1$
\sn
    \item  $\ortp_\gs(a_t,M_1,M_2)$ does not fork over $M_0$.
\sn
    \item  There is a sequence $\LL t(\alpha) : \alpha < \alpha(*)\RR$ listing $I$ with no repetitions, and a  $\le_\gs$-increasing sequence 
    $\LL M_{1,\alpha} : \alpha \le \alpha(*)+1\RR$ with $M_1 \le_\gs M_{1,0}$ and $M_2 \le M_{1,\alpha(*)+1}$ such that $a_{t(\alpha)} \in M_{1,\alpha +1}$ and
    $\ortp_\gs(a_{t(\alpha)},M_{1,\alpha},M_{1,\alpha +1})$ does not fork over $M_0$.
\end{enumerate}
\mn
4) Let $(M,N,\bfJ) \in K^{3,\bs}_\gs$ \underline{if} $M
\le_\gs N$ and $\bfJ$ is independent in $(M,N)$.

\noindent
5) Let $(M,N,\bfJ) \in K^{3,\qr}_\gs$ \underline{if}:
\mn
\begin{enumerate}
    \item  $M \le_\gs N$
\sn
    \item  $\bfJ$ is independent in $(M,N)$.
\sn
    \item If $M \le_\gs N'$ and $h$ is a one-to-one function from $\bfJ$ into $N'$ such that $(M,N',h''(\bfJ)) \in K^{3,\bs}_\gs$ \underline{then} there is a $\le_\gs$-embedding $g$ of $N$ into $N'$ over $M$ extending $h$.
\end{enumerate}
\end{definition}


\begin{remark} \label{g32}
We now can imitate relations of the axioms (as in
\cite[\S2]{Sh:600}), and basic properties of the notions introduced in
\ref{g29}. 
\end{remark}

\begin{definition}\label{g35}
1) We say $p$ is \emph{strongly dominated} by $\{p_t : t \in I\}$ and write 
$p \leq_\st \{p_t : t \in I\}$ (this set may
contain repetitions\footnote{So pedantically, we should use a sequence and 
    write $p \le_\st \LL p_t:t \in I\RR$.}) 
\underline{when}:
\mn
\begin{enumerate}  
    \item  $p \in \clS^\bs_\gs(N)$, $p_t \in \clS^\bs_\gs (N_t)$, 
    $N_t \le_\gs N^+ \in K_\gs$, $N \le_\gs N^+$ and
\sn
    \item  If $N^+ \le_\gs N^*$, $a_t \in N^*$,
    $\ortp(a_t,N^+,N^*) \in \clS^\bs_\gs(N^+)$ is parallel to $p_t$ 
    and $p' \in \clS^\bs_\gs(N^+)$ is parallel to $p$ (see Definition \ref{g11}), and $\{a_t : t \in I\}$ is independent in $(N^+,N^*)$ \underline{then} some $a \in N^*$ realizes $p'$.
\end{enumerate}
\mn
2) We say $p$ is weakly dominated by $\{p_t : t \in I\}$ and write 
$p \le_{\wk} \{p_t : t \in I\}$ \underline{when} for some set $J$
and function $h$ from $J$ onto $I$ we have $p \le_\st \{p_{h(t)} : t \in J\}$.

\sn
3) Let `dominated' mean strongly dominated.

\sn
4) We say $\gs$ is a strongly good $\iota$-frame \underline{when}
\begin{enumerate}
    \item It is a good $\iota$-frame.

    \item If $\bfJ$ is the disjoint union of $\bfJ_1$ and $\bfJ_2$, 
    $(M,N,\bfJ) \in K_\gs^{3,\bs}$, $M \leq_\gs M_1 \leq_\gs N$, and $(M,M_1,\bfJ) \in K_\gs^{3,\qr}$ \underline{then} $(M_1,N,\bfJ_2) \in K_\gs^{3,\bs}$ and $\ortp(a,M_1,N)$ does not fork over $M$ for all $a \in \bfJ_2$.
\end{enumerate}
\end{definition}

\begin{claim}\label{g38}
Assume $\gs$ is strongly good.

\sn
1) If $p$ is strongly dominated by $\{p_t : t \in I\}$ \underline{then} 
$p$ is weakly dominated by\\ $\{p_t : t \in I\}$.

\sn
2) If $p$ is strongly dominated by $\{p_t : t \in I\}$ \underline{then} 
for some $J \subseteq I$ of cardinality $< \kappa_\gs$, $p$ is strongly
dominated by $\{p_t : t \in J\}$. 

\sn
3) $p$ is weakly dominated by $\{p_t : t \in I\}$ \underline{iff} 
for some $\LL i_t : t \in I\RR$, $p$ is strongly dominated by 
$\big\{p'_s : s \in \{(t,i) : t \in I,\ i < i_t\}\big\}$, where 
$p'_{(t,i)} = p_t$ and $i_t < \kappa_\gs$ for each $t \in I$.

\sn
4) In Definition \ref{g35}(2) without loss of generality $(\forall s \in I)
(\exists^{< \kappa} t \in J)[h(t)=s]$. 

\sn
5) Preservation by parallelism.

\sn
6) $(M,N,\bfJ) \in K_\gs^{3,\bs}$ \underline{iff} for every finite 
$\bfI \subseteq \bfJ$ we have $(M,N,\bfI) \in K_\gs^{3,\bs}$.

7) If $\bfJ  = \{a_s : s \in I\}$ is independent in $(M,N)$ as witnessed by \\
$\LL M_{1,\alpha} : \alpha \leq \alpha(*)\RR$, $\LL t(\alpha) : \alpha < \alpha(*)\RR$ (see \ref{g29}(5)) then $(M_{1,0},M_{1,\alpha(*)},\bfJ) \in  K_\gs^{3,\qr}$.
\end{claim}

\begin{PROOF}{\ref{g38}}
1) Easy.

\sn
2) By \ref{g38}(2), following some manipulation. 

\sn
3) By \ref{g39} and clause (i) of Definition \ref{g2}.

\sn
4),5) Easy.

\sn
6) The $\Rightarrow$ direction follows from the definition; the $\Leftarrow$ direction can be proved by induction on $\alpha(*)$.

\sn
7) By Definition \ref{g35}(4).
\end{PROOF}

\begin{claim}\label{g39}
If $(M,N,\bfJ) \in K_\gs^{3,\bs}$ and $b \in N$ then there exist 
$\bfI \subseteq \bfJ$ and $M_1$ such that:
\begin{enumerate}
    \item $M \leq_\gs M_1 \leq_\gs N$

    \item $|\bfI| < \kappa_\gs$ 

    \item $b \in M_1$

    \item $(M,M_1,\bfI) \in K_\gs^{3,\bs}$
\end{enumerate}
\end{claim}

\begin{PROOF}{\ref{g39}}
Without loss of generality, $b \notin \bfJ$. We try, by induction on 
$i \leq \kappa$, to choose $N_i$ (and if possible, $\bfI_i$) such that
\begin{enumerate}
    \item[$(*)$]
    \begin{enumerate}
        \item $N_i \leq_\gs N$ and $\bfI_i \subseteq \bfJ \setminus \bigcup\limits_{j<i} \bfI_j$, with $|\bfI_i| < \kappa$.
    
        \item If $j < i$ then $N_j \leq_\gs N_i$ and $(N_j,N_{j+1},\bfI_j) \in K_\gs^{3,\pr}$.

        \item $N_0 = M$

        \item If $i$ is a limit ordinal then $N_i$ is prime over 
        $\LL N_j : j < i\RR$.

        \item If $i = j+1$ and $N_j$ has already been defined with $b \notin N_j$, and there is $\bfI\ \subseteq \bfJ \setminus \bigcup\limits_{\ell <j} \bfI_\ell$ of cardinality $< \kappa$ (or simply finite) such that 
        $$\big(N_j,N, \bfI \cup \{b\}\big) \notin K_\gs^{3,\bs}$$
        
        \underline{then} we can choose such $\bfI$ as our $\bfI_j$ and choose 
        $N_i \leq_\gs N$ such that $(N_j,N_i,\bfI_j) \in K_\gs^{3,\pr}$.
    \end{enumerate}
\end{enumerate}
If we carry the induction for all $i < \kappa$ we get a contradiction ({see} \ref{g2}(E)(c)), so for some $i(*) < \kappa$ we will hit a point where $N_{i(*)}$ is well defined, but $\bfI_{i(*)}$ is not.

We prove, by induction on $\theta \leq |\bfJ|$, that if $\bfI \subseteq \bfJ' = \bfJ \setminus \bigcup_{j<i(*)} \bfI_j$ has cardinality $\theta$ then $(N_{i(*)},N,\bfI \cup \{b\}) \in K_\gs^{3,\bs}$. So, using Case 1 of Definition \ref{g2}(E)(i), we are finished.
\end{PROOF}


\begin{claim}\label{g41}
1) If $p \le_{\wk} \{p_i:i < i^*\}$ and $i < i^* 
\Rightarrow q \perp p_i$ then $q \perp p$ (see Definition \ref{g26}(3)). 

\sn
2) If $p \le_{\wk} \{p_i:i < i^*\}$ and
$p \in \clS^\bs_\gs(M)$ \underline{then} $p \not\perp p_i$ for some $i < i^*$.

\sn
3) If $p \le_\st \{p_i : i < \alpha\}$
then $p \le_\st \{p_i : i < \alpha,\ p_i \not\perp p\}$
(see Definition \ref{g35}). 
\end{claim}

\begin{PROOF}{\ref{g41}}

1) By induction on $i^*$: for $i^*$ limit we use \ref{g2}(E)(i), and for $i^*$ successor use $q \perp p_{i^*-1}$.

\sn
2) By part (1) and \ref{g38}(3).

\sn
3) Easy.
\end{PROOF}

\begin{claim}\label{g44}
Assume $\gs$ is type-full.

If $\chi = \chi^{< \kappa} \in [\lambda,\mu)$, the following is impossible:
\begin{enumerate}
    \item[$(a)$]    $\LL M_i : i < \chi^+ \RR$ is $\le_\gs$-increasing $\le_\gs$-semi-continuous,
\sn
    \item[$(b)$]   $\LL N_i:i < \chi^+ \RR$ is $\le_\gs$-increasing, $\le_\gs$-semi-continuous,
\sn
    \item[$(c)$]   $M_i \le_\gs N_i \in K_{\le \chi}$,
\sn
    \item[$(d)$]  for some stationary 
    $S \subseteq \{\delta < \chi^+ : \cf(\gamma) \ge \kappa\}$, for every $i \in S$,
    \begin{itemize}
        \item  There is $a_i \in M_{i+1} \setminus M_i$ such that 
        $\ortp(a_i,N_i,N_{i+1})$ is not the non-forking extension of $\ortp(a_i,M_i,M_{i+1}) \in \clS^\bs_\gs(M_i)$.
    \end{itemize}
\end{enumerate}
\end{claim}

\begin{PROOF}{\ref{g44}}  
For some club $E$ of $\chi^+$, we have $i \in E \wedge j \in[i,\chi^+) \Rightarrow N_i \cap M_j = M_i$. For each $i \in S \cap E$, by \ref{g2}(E)(c),
there is a $j_i < i$ such that $\ortp(a_i,M_i,M_{i+1})$ does not fork over 
$M_{j_i}$. By clause (E)(i) of \ref{g2}, for some $j \in [j_i,i)$, we 
have that $\ortp(a_i,N_j,N_{i+1})$ is not the non-forking extension of $\ortp(a_i,M_{j_i},M_{i+1})$, so without loss of generality this holds for 
$j = j_i$. 

By Fodor's Lemma, for some $j(*) < j$ the set 
$S' = \{i \in S \cap E : j_i = j(*)\}$ is stationary. So 
$\{{b}_i : i \in S'\}$ is independent in $\big( \bigcup\limits_j M_j,M_{j+1} \big)$. By part (3) we are done.

Also, there is a sequence $\LL M_{j(*),\eps} : \eps \leq \eps_* \leq \kappa\RR$ which is $\eps_\gs$-increasing continuous, with $M_{j(*),0} = M_{j(*)}$, $M_{j(*),\eps} = N_{j(*)}$, and $(M_{j(*),\eps}, M_{j(*),\eps+1},c_\eps) \in K_\gs^{3,\pr}$. Now we can choose $\zeta_\eps < \chi^+$ by induction on 
$\eps < \eps*$, increasing continuous, such that 
$\big\{a_i : i \in [\zeta_i,\chi^+)\big\}$ is independent in $(M_{j(*),\eps},\bigcup\limits_j N_j)$ and $\ortp(a_i,M_{j(*),\eps},N_{i+1})$ does not fork over $M_{j(*)}$ for $i \in [\zeta,\chi^+)$ --- an easy contradiction. The induction works for $\eps = 0$ trivially, for $\eps$ limit by \ref{g38}(6), and for $\eps = \xi+1$ we use \ref{g39}.
\end{PROOF}

\begin{example}\label{g47}
1) For a complete f.o. strictly stable $T$ and $\kappa = \cf(\kappa) > \aleph_0$, we define $\gk$ by
\begin{enumerate}
    \item $K_\gk$ is the class of $\kappa$-saturated (equivalently, $\bfF_\kappa^a$-saturated) models of $T$.

    \item $\leq_\gk$ is defined by $M \leq_\gk N$ iff $M,N \in K_\gk$ 
    and $M \prec N$.
\end{enumerate}
\end{example}

\begin{claim}\label{g50}
If $p,p_i \in \clS^\bs_\gs(M)$ for $i < \kappa_\gs$ and 
$i < j \Rightarrow p_i \perp p_j$ \underline{then} 
$p \perp p_i$ for every $i < \kappa$ large enough.
\end{claim}

\begin{PROOF}{\ref{g50}}
 Similar to the proof of \cite[1.6=Lj20]{Sh:1239}.   
\end{PROOF}

\begin{definition}\label{g53}
1) We say that a good frame $\gs$ is $\theta$-based$_1$ \underline{when}:
\mn
\begin{enumerate}
    \item   If $p \in \clS^\bs_\gs(M)$ then for some 
    $\bar\bfa \in {}^{\theta >}\!M$, $p$ is based on $\bar\bfa$
    (see Definition \ref{g11}(4)).
\end{enumerate}
\mn
2) We say that $\gs$ is $\theta$-based$_2$ \underline{when}:
\mn
\begin{enumerate}
    \item  Is as in part (1).
\sn
    \item  $\gs$ is type-full.
\sn
    \item  If $M_1 \le_\gs M_2$ and $p \in \clS^\bs_\gs(M_2)$ \underline{then}, for some $\bar a_\ell \in {}^{\theta >}(M_\ell)$, 
    the types $p$ and $\ortp_\gs(\bar a_2,M_1,M_2)$ are based on $\bar\bfa_2,\bar\bfa_1$, respectively.
\end{enumerate}
\end{definition}
\newpage

\section{Thoughts on the main gap} \label{7}

We address here two problems: type theory (i.e. dimension, orthogonality,
etc.) for strictly stable classes, and the main gap concerning somewhat
saturated models. The hope always was that advances in the first will
help the second.

Concerning the first order case, work started in \cite[Ch.V]{Sh:c}
(particularly \S5) and \cite{Sh:429} and was much
advanced in Hernandes \cite{He92}; but this was not enough for the main
gap for somewhat saturated models.

Here we are dealing with the type dimension in a general framework.
\bigskip

\centerline {$* \qquad * \qquad *$}
\bigskip

The main gap for $\aleph_1$-saturated models of a countable first
order theory is open.  \emph{A priori}, it has looked easier than the one for
models (which was preferred, being ``the original question") because of
the existence of prime models over any, but is still
open. (The problem for uncountable first order $|T|^+$-saturated models is as well).

Why doesn't the proof in \cite[Ch.XII]{Sh:c} work?
What's missing is, in $\gC^{\eq}$,
\mn
\begin{enumerate}
    \item[$\circledast$]   If $M_0 \prec  M_1 \prec M_2$ are $\aleph_1$-saturated, $a \in M_2 \setminus M_1$ and $(a/M_1) \not\perp M_0$ \underline{then} for some $b \in M_2 \setminus M_1$ we have $\nonfork{b}{M_1}_{M_0}$.
\end{enumerate}
\mn
The central case is when $a/M_1$ is orthogonal to $q$ if $q \perp M_0$.

\bn
\textbf{Possible Approach 1}:  We use $T$ being first order
countable, stable NDOP (even shallow) to understand types.  
See \cite{Sh:851}.

\bn
\textbf{Possible Approach 2}:  We use the context {dealt with in this paper}. 
We are poorer in knowledge on the class \underline{but} we have a richer 
$\gC^{\eq}$, so we may prove $\circledast$ even if it fails for $T$ in the
elementary case (this is a connection between \cite{Sh:839} and this work).

\bn
\textbf{Possible Approach 3}:  We start with the context {here}.  If things
are not OK, we define such a derived DAEC; this was done in 
\cite{Sh:300f} and \cite{Sh:600}. It may have non-structure properties --- 
enough to get the maximal number of models up to isomorphism. If not, we 
arrive to a finer $\gk$, but still a case of our context.  Similarly in limit.
If we succeed enough times we shall prove that all is OK.

\bn
\textbf{Possible Approach 4}:  Now we have a maximal non-forking tree
$\LL M_\eta,a_\eta:\eta \in \clT\RR$ inside a somewhat
saturated model; for \cite{Sh:c}, e.g. $\|M_\eta\| \le \lambda$, the models
are $\lambda^+$-saturated but we use models from 
here.  If $M$ is prime over $\bigcup\{M_\eta : \eta \in \clT\}$ we are done, 
but maybe there is a residue. This appears in the following way: 
for $\eta \in \clT$ and $p \in \mathbf \clS^\bs(M_\eta)$, 
the dimension of $p$ is not exhausted by 
$$\{ a_{\eta \caret\LL \alpha \RR}:\eta \caret \LL\alpha \RR \in \clT 
\text{ and } (a_{\eta \caret \LL \alpha \RR}/M_\eta) \not\perp p\}$$
but the lost part is not infinite!  This imposes $\le \lambda$ unary
functions from $\clT$ to $\clT$.  Now it seems to us that the
question of {whether} this possible non-exhaustion {can}   
    arise\footnote{Essentially: there is a non-algebraic $p \in (M^\bot)^\bot$ which do not 1-dominate any $q \in \clS(M)$.} 
is not a good dividing line, as though its negation is
informative it is not clear whether it has any consequence.  However,
there are two candidates for dividing lines (actually, their
disjunction seems {to be what we want}).
\mn
\begin{enumerate}
    \item[$(A)$]   $(*) \quad$ We can find $M$, 
    $\LL M_\eta,a_\eta : \eta \in \clT\RR$ as above and 
    $\eta_* \in \clT$, $\lh(\eta) = 2$, $\nu_* \in \clT$, $\lh(\nu_*) = 1$, $\eta_* \rest 1 \ne \nu_*$, and $p \in \clS^\bs(M_{\eta_*})$, 
    $p \perp M_{\eta \rest 1}$ with a residue as above such that we need $M_{\nu_*}$ to explicate it.
\end{enumerate}
\mn
More explicitly,
\mn
\begin{enumerate}
    \item[$(*)'$]   If $M' \le_\gs M$ is prime over 
    $\bigcup\{M_\eta : \eta \in \clT\}$ and we can find 
    $a_{\eta_*,\nu_*} \in M \setminus M'$ such that $\ortp\big(\cC(a_{\eta^*,\nu^*},M'),\bigcup\{M_\eta:\eta \in \clT\}\big)$ mark $(M_{\eta_*},M_{\nu_*})$.
\end{enumerate}
\mn
Even in $(*)'$ we have to say more in order to succeed in using it.

From $(*)'$ we can prove a non-structure result: on $\clT$ we can
code any two-place relation $R$ on 
$\{\eta \in \clT : \lh(\eta) = 1,\ M_\eta,M_{\eta_* \rest 1} 
\text{ isomorphic over } M_{\LL\ \RR}\}$ which is {of the form} 
$\nu_1\ R\ \nu_2 \Leftrightarrow (\exists \nu) \bigwedge \limits_{\ell}$
[there is $\eta',\ \eta_\ell \triangleleft \eta' \in \clT,\ \lh(\eta') = 2$ and $\nu \in T,\ \lh(\nu)=1$ and there is $a_{\eta',\nu}$ as above].

More complicated is the case
\mn
\begin{enumerate}
    \item[$(B)$]  $(**) \quad$ We can fix 
    $M$, $\LL M_\eta,a_\eta : \eta \in \clT\RR$ as above, 
    $\eta_* \in T$, $\nu_*,\nu_{**} \in \clT$, $\lh(\eta_*) = \lh(\nu_*) = \lh(\nu_{**}) = 1$ such that $(\eta_*,\nu_*),(\eta_*,\nu_{**})$ are as above.
\end{enumerate}
\mn
But whereas for (A) we have to make both $\eta_*$ and $\nu_*$ not
redundant in (B), in order to get non-structure we have to use a case of (B)
which is not ``a faking;" e.g. we cannot replace
$(M_{\eta_*},a_{\eta_*})$ by two such pairs.

That is, the ``faker" is a case where we can find
$M'_{\eta_*},M''_{\eta_*}$ such that:
\mn
\begin{itemize}
    \item   $\NF(M_{\LL\ \RR},M'_{\eta_*},M''_{\eta_*},M_{\eta_*})$
\sn
    \item  $M_{\eta_*}$ is prime over $M'_{\eta_*} \cup M''_{\eta_*}$.
\sn
    \item  Only $(M'_{\eta_*},M_{\nu_*})$ and $(M''_{\eta_*},M_{\nu_{**}})$ relate.
\end{itemize}

\begin{enumerate}
    \item[(C)] If both (A) and (B), in the right formulation, do not appear then
    \begin{enumerate}
        \item[$(\alpha)$]  \textbf{A good possibility} 
        
        We can prove a structure theory: for $M$, $\LL M_\eta,a_\eta:\eta \in \clT\RR$ as above; that is, on each $\suc_\clT(\eta)$ we have a two-place relation, but it is very simple: you have to glue some together or {expand the set of successors by a tree structure}.
    \end{enumerate}
\end{enumerate}
If this fails we may fall back to approach (3).

We may consider (see \cite{Sh:897}, \cite{Sh:1133}):

\begin{question}  \label{z3}
1) For an AEC $\gk$, when does the theory of a model in the logic ${\cL} =
\bbL_{\infty,\kappa} [\gk]$ enriched by dimension quantifiers, 
characterize models of $\gk$ up to isomorphism?  Similarly
enriching also by game quantifiers of length $\le \kappa$.

\sn
2) Prove the main gap theorem in the version: if $\gs$ is
$n$-beautiful \saharon{or $n+1$?} \underline{then} for 
$K_{\lambda^{+n}}$ the main gap holds. 
In particular, if $\gs$ has NDOP, then every $M \in K_{\lambda^{+n}}$ 
is prime over some non-forking tree of 
$\le_{{\gK}[{\gs}]}$-submodels $\LL M_\eta:\eta \in \clT\RR$,
each $M_\eta$ of cardinality $\le \lambda$, where 
$\clT \subseteq {}^{\omega >}(\lambda^{+n})$.  If $\gs$ is shallow then the tree
has depth $\le \Depth({\gs}) < \lambda^+$ and we can draw a
conclusion on the number of models.
\end{question}

\begin{discussion}\label{h5}
Assume stability in $\lambda_\gs$.

Let $M_0 \in K_\gs$, $\lambda^+_\gs$-saturated, at least for the time
being.

\noindent
1) Assume
\mn
\begin{enumerate}
    \item[$\boxplus_1$]  $N_0 \le_\gs N_1 \le_\gs M$, $N_\ell \in K^\gs_\lambda$, $a \in N_0$, and $(N_0,N_1,a) \in K^{3,\pr}_\gs$.
\end{enumerate}
\mn
We choose $(N^+_{1,i},N_{1,i},\bfI_i)$ and also, if possible, 
$(M_1,a_i)$ by induction on $i \le \lambda^+_\gs$ such that
\mn
\begin{enumerate}
    \item[$(*)$]  
    \begin{enumerate}
        \item  $N_{0,i} \le_\gs N_{1,i} \le_\gs N^+_{1,i} \le_\gs M$

        \item  $\bfI_i \subseteq \{c \in M : \ortp(c,N_{1,i},M_0) \perp N_0\}$ is independent in $(N_{1,i},N^+_{1,i},M)$ and minimal.

        \item  $\LL N_j:j \le i\RR$ is $\le_\gs$-semi-continuous; also, $\LL N^+_j:j \le i\RR$ {is as well.}

        \item If $i=j+1$ then $N^+_{1,i}$ is $\le_\gs$-universal over $N^+_{1,j}$ and $(N_0,N_{1,i},a) \in K^{3,\pr}_\gs$.

        \item If $j < i$ then $\bfI_j \setminus (N_i \cap \bfI_j) \subseteq \bfI_i$.

        \item If possible:
        \begin{enumerate}
            \item[$(\alpha)$] $N_i \le_\gs M^+_i \le_\gs M$

            \item[$(\beta)$]  $(\bfI_i \setminus M_i)$ is independent in $(M_i,M)$.

            \item[$(\gamma)$] $a_i \in M \setminus (\bfI_i)$

            \item[$(\delta)$] $\ortp(a_i,M^*_1,M) \in \clS^\bs_\gs(N^+_i)$ is $\perp N_i$.

            \item[$(\eps)$] $N^*_i \le N_{1,i+1}$  
        \end{enumerate}

        \item If $i=j+1$ {and} there are $(b,N^+_*,N_{**})$ such that 
        $b \in N^+_{1,j} \setminus N_{1,j}$, 
        $$N_{1,i} \le_\gs N_* \le_\gs N_{**} \in K^\gs_{\lambda_\gs},$$  $N^+_{1,i} \le_\gs N_{**}$, and $\ortp_\gs(b,N_*,N_{**})$ forks over $N_{1,j}$ \underline{then}, for some 
        $b \in N^+_{1,j} \setminus N_{1,j}$, the type $\ortp_\gs(N_{1,i},N^+_{1,i})$ forks over $N_{1,j}$.
    \end{enumerate}
\end{enumerate}
\mn
There is no problem to carry the induction.
\mn
\begin{enumerate}
    \item[$\boxplus_2$]  The following subset of $\lambda^+_\gs$ is not  stationary --- say, disjoint to the club $C$:
    \begin{itemize}
        \item  $S = \{i < \lambda^+_\gs:\cf(i) \ge \kappa_\gs$ and $(M_i,a_i)$ is well defined$\}$
\sn
        \item $S_2 = \big\{i : \cf(i) \ge \kappa_\gs$ and for some $b \in N^+_{1,i},\ \tp(b,N_{1,i},N^+_{1,i}) = N_0 \big\}$.
    \end{itemize}
\end{enumerate}
\mn
2) Similarly without $(N_0,a)$ hence without ``$\perp N_0$;'' {it's} just simpler.
\end{discussion}

\begin{definition}\label{h8}
We say $(\olsi N,\bar a,\bar I)$ is a decreasing pair for $M$ \underline{when}
for some $n$:
\mn
\begin{enumerate}
    \item  $\olsi N = \LL N_\ell : \ell \le n \RR$ is $\le_\gs$-increasing. 
\sn
    \item  $N_\ell \le_\gs M$, $N_\ell \in K^\gs_{\lambda_\gs}$
\sn
    \item  $\bar a = \LL a_\ell:\ell < n\RR$
\sn
    \item  $(N_\ell,N_{i+1},a_\ell) \in K^{3,\pr}_\gs$
\sn
    \item  $\bar\bfI = \LL \bfI_\ell:\ell \le n \RR$
\sn
    \item  $\bfI_\ell$ is independent in $(N_\ell,M)$.
\sn
    \item  $\bfI_\ell \subseteq \{c \in M:\ortp(c,N_\ell,M) \in \clS^\bs_\gs(N_\ell)$ is $\perp N_k$ if $k < \ell\}$
\sn
    \item  If $N_\ell \le_\gs N \le_\gs M$, $b \in M \setminus N_0 \setminus \bfI_\ell$, and $\ortp(b,N,M)$ is $\not\perp N_\ell$ but is orthogonal to $N_k$ for $k < \ell$ \underline{then} $b$ depends on $\bfI_\ell$ in $(N_\ell,M)$.
\end{enumerate}
\end{definition}

\noindent
\textbf{Attempt to prove decomposition}

We assume dimensional continuity to prove decomposition.  If we {would} like
to get rid of ``$M$ is $\lambda^+_\gs$-saturated", we {must} assume we have
a somewhat weaker version $\gs_*$ of $\gs$ where $\lambda_{\gs_*} <
\lambda_\gs$ and $\LL N_{0,i} : i < \lambda_\gs\RR$ $\le_{\gs_*}\!$-represent $N_0$, and work with that.  Assuming CH, $|T| =
\aleph_0$ is fine.  Without dimensional discontinuity we call `nice' {any} $(\bar N,\bar a,\bar\bfI)$ of length $\le \kappa_\gs$!
\bigskip

\centerline {$* \qquad * \qquad *$}
\bigskip

\begin{definition}\label{h11}
We say $\bfd = (I,N,\bar a,\bar\bfI) = (I_\bfd,\olsi N_d,\bar a_\bfd,\bar\bfI_\bfd)$ is a partial
  decomposition of \underline{when}:
\mn
\begin{enumerate}
    \item[$\boxplus$]
    \begin{enumerate}
        \item $I \subseteq {}^{\omega >}\Ord$ is closed under initial segments.

        \item  $\olsi N = \LL N_\eta:\eta \in I\RR$, so $N_\eta = N_{\bfd,\eta}$.

        \item  $\bar a = \big\LL a_\eta : \eta \in I\setminus \{\LL\ \RR\} \big\RR$, so $a_\eta = a_{\bfd,\eta}$.

        \item $\bar\bfI = \LL \bfI_\eta : \eta \in I\RR$, so $\bfI_\eta = \bfI_{\bfd,\eta}$. 

        \item If $\eta \in I$ then $$\big( \big\LL N_{\eta \rest \ell} : \ell \le \lh(\eta) \big\RR, \big\LL \bar a_{\eta \rest (\ell +1)} : \ell < \lh(\eta) \big\RR, \big\LL \bfI_{\eta \rest \ell} : \ell \le \lh(\eta)\big\RR \big)$$ is nice in $M$.

        \item If $\eta \in I$ then $\LL a_{\eta \caret \LL \alpha \RR}:\eta \caret \LL \alpha \RR \in I\RR$ is a sequence of members of $\bfI_\eta$ with no repetitions.
    \end{enumerate}
\end{enumerate}
\end{definition}

\begin{definition}\label{h14}
Let $\le_\mu$ be the following two-place relation on the set of
decompositions of $M$:

$\bar\bfd_1 \le_M \bfd_2$ \underline{iff}
\mn
\begin{enumerate}
    \item $I_{\bfd_1} \subseteq \bfI_{d_1}$
\sn
    \item  $\olsi N_{\bfd_1} = \olsi N_{\bfd_2} \rest I_{d_1}$
\sn
    \item  $\bar a_{\bfd_1} = \bar a_{\bfd_2} \rest (I_{\bfd_1} \setminus \{<j\})$
\sn
    \item  $\bar\bfI_{\bfd_1} = \bar\bfI_{d_2} \rest I_{\bfd_1}$ 
\end{enumerate}
\end{definition}

\bibliographystyle{amsalpha}
\bibliography{shlhetal}
\end{document}

%% file: mgrimes.tex
\newcommand{\Car}{\mathrm{Car}}

\newcommand{\comp}{\mathrm{comp}}
\newcommand{\dist}{\mathrm{dist}}

\newcommand{\nmet}{\mathrm{nmet}}
\newcommand{\na}{\mathrm{na}}

\newcommand{\NF}{\mathrm{NF}}
\newcommand{\ortp}{\mathrm{ortp}}

\newcommand{\qr}{\mathrm{qr}}

\newcommand{\rnd}{\mathrm{rnd}}
\newcommand{\st}{\mathrm{st}}

\newcommand{\wk}{\mathrm{wk}}

\newcommand{\cf}{\mathrm{cf}}
\newcommand{\cl}{{c\ell}}

\newcommand{\Depth}{\mathrm{Depth}}

\newcommand{\id}{\mathrm{id}}

\newcommand{\Ord}{\mathrm{Ord}}

\newcommand{\at}{\mathrm{at}}

\newcommand{\bs}{\mathrm{bs}}

\newcommand{\EC}{\mathrm{EC}}

\newcommand{\eq}{\mathrm{eq}}

\newcommand{\LST}{\mathrm{LST}}

\newcommand{\pr}{\mathrm{pr}}

\newcommand{\suc}{\mathrm{suc}}

\newcommand{\tp}{\mathrm{tp}}

\newcommand{\bfa}{\mathbf{a}}

\newcommand{\bfd}{\mathbf{d}}

\newcommand{\bfF}{\mathbf{F}}

\newcommand{\bfI}{\mathbf{I}}

\newcommand{\bfJ}{\mathbf{J}}

\newcommand{\bbL}{\mathbb{L}}

\newcommand{\bbQ}{\mathbb{Q}}

\newcommand{\mn}{\medskip\noindent}
\newcommand{\sn}{\smallskip\noindent}
\newcommand{\bn}{\bigskip\noindent}

\newcommand{\cC}{\mathscr{C}}

\newcommand{\cE}{\mathscr{E}}

\newcommand{\cL}{\mathscr{L}}

\newcommand{\cU}{\mathscr{U}}

\newcommand{\clP}{\mathcal{P}}

\newcommand{\clS}{\mathcal{S}}
\newcommand{\clT}{\mathcal{T}}

\newcommand{\gC}{\mathfrak{C}}

\newcommand{\gK}{\mathfrak{K}}

\newcommand{\gk}{\mathfrak{k}}

\newcommand{\gs}{\mathfrak{s}}

\newcommand{\eps}{\varepsilon}
\newcommand{\lh}{{\ell g}}
\newcommand{\rest}{\restriction}

\newcommand{\caret}{{\char 94}}
\newcommand{\LL}{\langle}
\newcommand{\RR}{\rangle}

\newcommand{\olsi}[1]{\,\overline{\!{#1}}} 

\newcommand*{\defeq}{\mathrel{\vcenter{\baselineskip0.5ex \lineskiplimit0pt\hbox{\scriptsize.}\hbox{\scriptsize.}}}=}

\newcommand{\saharon}[1]{\magenta{\textbf{[#1]}}}

\newcount\skewfactor
\def\mathunderaccent#1#2 {\let\theaccent#1\skewfactor#2
\mathpalette\putaccentunder}
\def\putaccentunder#1#2{\oalign{$#1#2$\crcr\hidewidth
\vbox to.2ex{\hbox{$#1\skew\skewfactor\theaccent{}$}\vss}\hidewidth}}

\newbox\noforkbox \newdimen\forklinewidth
\setbox0\hbox{$\textstyle\bigcup$}
\setbox1\hbox to \wd0{\hfil\vrule width \forklinewidth depth \dp0
                        height \ht0 \hfil}
\setbox\noforkbox\hbox{\box1\box0\relax}
\def\unionstick{\mathop{\copy\noforkbox}\limits}
\def\nonfork#1#2_#3{#1\unionstick_{\textstyle #3}#2}
\def\nonforkin#1#2_#3^#4{#1\unionstick_{\textstyle #3}^{\textstyle
    #4}#2}
\setbox0\hbox{$\textstyle\bigcup$}
\setbox1\hbox to \wd0{\hfil{\sl /\/}\hfil}
\setbox2\hbox to \wd0{\hfil\vrule height \ht0 depth \dp0 width
                                \forklinewidth\hfil}
\newbox\doesforkbox
\setbox\doesforkbox\hbox{\box1\box0\relax}
\def\nunionstick{\mathop{\copy\doesforkbox}\limits}

\def\fork#1#2_#3{#1\nunionstick_{\textstyle #3}#2}
\def\forkin#1#2_#3^#4{#1\nunionstick_{\textstyle #3}^{\textstyle
    #4}#2}

\newcommand{\stickT}{%
\setbox255=\hbox{\raise1ex\hbox{$\hspace{0.2pt}\,\bullet\,$}}
\mathord{\rlap{\hbox to\wd255{\hss\hbox{$|$}\hss}}
\box255}
}
\newcommand{\stickS}{%
\setbox255=\hbox{\raise0.6ex\hbox{$\scriptstyle\bullet$}}
\mathord{\rlap{\hbox to\wd255{\hss\hbox{$\scriptstyle|$}\hss}}
\box255}
}